\documentclass[reqno,10pt,centertags,draft]{amsart}
\usepackage{amsmath,amsthm,amscd,amssymb,eufrak,latexsym,upref}

\makeatletter
\def\theequation{\@arabic\c@equation}

\newcommand{\diag}{\operatorname{diag}}

\newcommand{\bbN}{{\mathbb{N}}}
\newcommand{\bbR}{{\mathbb{R}}}

\newcommand{\bbC}{{\mathbb{C}}}

\newcommand{\cA}{{\mathcal A}}
\newcommand{\cB}{{\mathcal B}}
\newcommand{\cC}{{\mathcal C}}

\newcommand{\cF}{{\mathcal F}}
\newcommand{\cG}{{\mathcal G}}
\newcommand{\cH}{{\mathcal H}}
\newcommand{\cI}{{\mathcal I}}
\newcommand{\cJ}{{\mathcal J}}

\newcommand{\cL}{{\mathcal L}}
\newcommand{\cM}{{\mathcal M}}
\newcommand{\cN}{{\mathcal N}}

\newcommand{\cP}{{\mathcal P}}
\newcommand{\cQ}{{\mathcal Q}}
\newcommand{\cR}{{\mathcal R}}

\newcommand{\cU}{{\mathcal U}}
\newcommand{\cV}{{\mathcal V}}

\newcommand{\cY}{{\mathcal Y}}
\newcommand{\cZ}{{\mathcal Z}}

\newcommand{\gh}{\mathfrak{h}}
\renewcommand{\gg}{\mathfrak{g}}
\newcommand{\gG}{\mathfrak{G}}
\newcommand{\gH}{\mathfrak{H}}
\newcommand{\gR}{\mathfrak{R}}

\newcommand{\no}{\nonumber}
\newcommand{\lb}{\label}
\newcommand{\f}{\frac}
\newcommand{\ul}{\underline}
\newcommand{\ol}{\overline}

\newcommand{\wti}{\widetilde  }
\newcommand{\Oh}{O}

\newcommand{\loc}{\text{\rm{loc}}}
\newcommand{\spec}{\text{\rm{spec}}}
\newcommand{\rank}{\text{\rm{rank}}}

\newcommand{\dom}{\text{\rm{dom}}}
\newcommand{\ess}{\text{\rm{ess}}}
\newcommand{\res}{\text{\rm{res}}}

\newcommand{\ac}{\text{\rm{ac}}}
\newcommand{\singc}{\text{\rm{sc}}}
\newcommand{\p}{\text{\rm{p}}}
\newcommand{\supp}{\text{\rm{supp}}}

\newcommand{\AC}{\text{\rm{AC}}}
\newcommand{\bi}{\bibitem}
\newcommand{\hatt}{\widehat}
\newcommand{\kdv}{\text{\rm{KdV}}}
\newcommand{\skdv}{\text{\rm{s-KdV}}}

\renewcommand{\Re}{\text{\rm Re}}
\renewcommand{\Im}{\text{\rm Im}}

\numberwithin{equation}{section}

\newtheorem{theorem}{Theorem}[section]
\newtheorem{lemma}[theorem]{Lemma}
\newtheorem{corollary}[theorem]{Corollary}
\newtheorem{hypothesis}[theorem]{Hypothesis}
\newtheorem{definition}[theorem]{Definition}
\newtheorem{remark}[theorem]{Remark}

\newcommand{\abs}[1]{\lvert#1\rvert}

\begin{document}

\title[Matrix Schr\"odinger Operators with Finite-Band Spectra]{A 
Class of Matrix-Valued Schr\"odinger Operators with Prescribed 
Finite-Band Spectra}
\author[F.\ Gesztesy and L.\ A.\ Sakhnovich]{Fritz Gesztesy and Lev
A.~Sakhnovich}
\address{Department of Mathematics,
University of Missouri, Columbia, MO 65211, USA}
\email{fritz@math.missouri.edu}
\urladdr{http://www.math.missouri.edu/people/fgesztesy.html}
\address{735 Crawford Ave., Brooklyn, NY 11223, USA}
\email{Lev.Sakhnovich@verizon.net}
\date{\today}
\subjclass{}
\keywords{Matrix-valued Schr\"odinger operators, finite-band spectra,
Weyl--Titchmarsh matrices.}

\begin{abstract}
We construct a class of matrix-valued Schr\"odinger operators with
prescribed finite-band spectra of maximum spectral multiplicity. The
corresponding matrix potentials are shown to be stationary solutions
of the KdV hierarchy. The methods employed in this paper rely on
matrix-valued Herglotz functions, Weyl--Titchmarsh theory, pencils of
matrices, and basic inverse spectral theory for matrix-valued
Schr\"odinger operators.
\end{abstract}

\maketitle

\section{Introduction} \lb{s1}

While basic aspects of inverse spectral theory for matrix-valued
Schr\"odinger operators were established some time ago, finer
properties such as isospectral sets (manifolds) of potentials, for
instance, in the periodic or algebro-geometric finite-band cases, are
still in their infancy. This paper intends to make a modest contribution
to this circle of ideas. More precisely, given a closed set
$\Sigma\subset\bbR$ of the type
\begin{equation}
\Sigma=\Bigg\{\bigcup_{j=0}^{n-1} [E_{2j},E_{2j+1}]\Bigg\}\cup
[E_{2n},\infty),  \lb{1.1}
\end{equation}
where 
\begin{equation}
\{E_\ell\}_{0\leq \ell\leq 2n}\subseteq \bbR, \; n\in\bbN, \text{ with
$E_\ell<E_{\ell+1}$, $0\leq \ell\leq 2n-1$,} \lb{1.2}
\end{equation}
we construct $m\times m$ matrix-valued Schr\"odinger operators
$H_\Sigma=-d^2/dx^2 \cI_m+\cQ_\Sigma$ in $L^2(\bbR)^{m\times m}$ ($\cI_m$
the identity matrix in $\bbC^{m\times m}$, $m\in\bbN$) with prescribed
spectrum $\Sigma$,
\begin{equation}
\spec(H_\Sigma)=\Sigma, \lb{1.3}
\end{equation}
of uniform spectral multiplicity $2m$. The constructed matrix potentials
$\cQ_\Sigma$ will turn out to be reflectionless in the sense discussed in
\cite{CGHL00}, \cite{GT00}, and \cite{KS88}, that is, the
half-line Weyl--Titchmarsh matrices $\cM_{\pm,\Sigma}(z,x)$ associated
with
$H_\Sigma$, the half-lines $[x,\pm\infty)$, and a Dirichlet boundary
condition at $x\in\bbR$, satisfy 
\begin{align}
&\lim_{\varepsilon\downarrow 0}\cM_{+,\Sigma}(\lambda+i\varepsilon,x)
=\lim_{\varepsilon \downarrow 0}\cM_{-,\Sigma}(\lambda-i\varepsilon,x),
\lb{1.4} \\
& \hspace*{7mm} \lambda\in \bigcup_{j=0}^{n-1} (E_{2j},E_{2j+1}) \cup
(E_{2n},\infty), \; x\in\bbR. \no
\end{align}
Especially, $\cM_{+,\Sigma}(\cdot,x)$ is the analytic continuation 
of $\cM_{-,\Sigma}(\cdot,x)$ through the set $\Sigma$, and vice versa.
In other words, $\cM_{+,\Sigma}(\cdot,x)$ and $\cM_{-,\Sigma}(\cdot,x)$
are the two branches of an analytic matrix-valued function
$\cM_{\Sigma}(\cdot,x)$ on the two-sheeted Riemann surface of
$\big(\prod_{\ell=0}^{2n} (z-E_\ell)\big)^{1/2}$. These facts imply a 
purely absolutely continuous spectrum $\Sigma$ of the associated
Schr\"odinger operator $H_\Sigma$ of uniform (maximal) multiplicity
$2m$.  

Before we turn to a brief description of the contents of each section, it
seems appropriate to mention some of the pertinent results and
especially, the most recent activities in connection with (inverse)
spectral theory of matrix-valued Schr\"odinger operators. The basic
Weyl--Titchmarsh theory of singular Hamiltonian systems and
their basic spectral theory were developed by Hinton and Shaw, Kogan
and Rofe-Beketov, Orlov, and others (see, e.g., \cite{Hi69},
\cite{HS81}--\cite{HS86}, \cite{JNO00}, \cite{KR74}, \cite{Kr89a},
\cite{Kr89b}, \cite{Or76}, \cite{Sa88a}, \cite{Sa92}, 
\cite{Sa97b}, \cite{Sa94a}, \cite[Ch.~9]{Sa99a}, and the
references therein). Various aspects of direct spectral theory, including 
investigations of the nature of the spectrum involved, (regularized)
trace  formulas, uniqueness theorems, etc., appeared in \cite{BL00},
\cite{Ca99}--\cite{Ca01}, \cite{Cl92}, \cite{CGHL00},
\cite{CH98}, \cite{GH97}, \cite{GKM01}, \cite{GS00}, 
\cite{Kr83}, \cite{Pa95}, \cite{RH77}, \cite{RK01}. General
asymptotic expansions of Weyl--Titchmarsh matrices as the  (complex)
spectral parameter tends to infinity under optimal regularity assumptions
on the coefficients are of relatively recent origin and can be found in
\cite{Cl92}, \cite{CG01} (see also \cite{Sa88}, \cite{Tr00}). The
inverse scattering formalism for Schr\"odinger operators has
been studied by a variety of authors and we refer, for instance, to
\cite{AM63}, \cite{AG95}, \cite{AG98}, \cite{AK90}, \cite{MO82},
\cite{NJ55}, \cite{Ol85}, \cite{WK74}. General inverse spectral theory,
the existence of transformation operators, etc., are discussed in
\cite{Ma94}, \cite{Ma99}, \cite{Ro60}, \cite{Sa97b}, \cite{Sa94a},
\cite{Sa96}, \cite{Sa99}, \cite{Sa99a}, and the references therein.
Inverse monodromy problems have recently been discussed in \cite{AD97},
\cite{AD00}, \cite{AD00a}, \cite{Ca01}, \cite{Ma94}, \cite{Ma99},
\cite{Sa94}, \cite{Sa99a}, and the literature cited therein. More
specific inverse spectral problems, such as compactness of the
isospectral set of periodic Schr\"odinger operators
\cite{Ca00},  special isospectral matrix-valued Schr\"odinger operators,
and  Borg-type uniqueness theorems (for periodic coefficients as well as
eigenvalue problems on compact intervals) were recently studied in
\cite{Ch99a}, \cite{CLW01}, \cite{CS97}, \cite{CG01}, \cite{CGHL00},
\cite{De95}, \cite{JL98}, \cite{JL99}, \cite{Ma94}, \cite{Ma99},
\cite{Sa99a}, \cite{Sh00}, \cite{Sh01}, \cite{SS98}. Moreover, direct
spectral theory in the particular case of periodic Schr\"odinger operators
(i.e., Floquet theory and the like) has been studied in \cite{Ca99},
\cite{Ca00a}, \cite{CGHL00}, \cite{De80}, \cite{De95}, \cite{GL87},
\cite{Jo87}, \cite{Kr83a}, \cite{Kr83}, \cite{Ro63}, \cite{Sh00},
\cite{Sh01}, \cite{Ya92}--\cite{YS75b}, with many more pertinent
references to be found therein. Apart from Floquet theoretic applications
in connection with Schr\"odinger operators already briefly touched upon,
we also need to mention applications to random Schr\"odinger operators
associated with strips as discussed, for instance, in \cite{CL90},
\cite{JNO00},
\cite{KS88}, and especially to nonabelian completely integrable systems.
Since the literature associated with the latter topic is of enormous
proportions, we can only refer to a few pertinent publications, such as,
\cite{AK90}, \cite{AS00}, \cite{CD77}, \cite{Di91}, \cite{Du83},
\cite{Ma78}, \cite{Ma88}, \cite{MO82}, \cite{OMG81}, \cite{OS98},
\cite{Po01}, \cite{Sa97a}, \cite{Sa88}--\cite{Sa94a}, \cite{Sa97}. The
interested reader will find a wealth of additional material in these
references.

Section \ref{s2} summarizes basic results in Weyl--Titchmarsh theory and
some elements of inverse spectral theory for matrix-valued Schr\"odinger
operators. Polynomial pencils of matrices are briefly reviewed in Section
\ref{s3}. In Section \ref{s4} we present our principal new result, the
construction of $m\times m$ matrix-valued Schr\"odinger operators
$H_\Sigma$ with spectrum $\Sigma$ (cf.\ \eqref{1.1}) and uniform maximal
spectral multiplicity $2m$. In our final Section \ref{s5} we prove that
$\cQ_\Sigma$ satisfies a stationary KdV equation (in fact, we explicitly
identify the first equation in the stationary KdV hierarchy satisfied by
$\cQ_\Sigma$) and derive matrix-valued trace formulas for $\cQ_\Sigma$ 
and higher-order KdV invariants. 

\section{Basic Facts on Weyl--Titchmarsh Theory} \lb{s2}

In this section we briefly recall basic elements of the Weyl--Titchmarsh 
theory for matrix-valued Schr\"{o}dinger operators. Throughout
this paper all matrices will be considered over the field of
complex numbers $\bbC$, and the corresponding linear space of
$k\times\ell$ matrices will be denoted by $\bbC^{k\times\ell}$,
$k,\ell\in\bbN$.  Moreover, $\cI_k$ denotes the identity matrix in
$\bbC^{k\times k}$, $\cM^*$ the adjoint (i.e., complex conjugate
transpose), $\cM^t$ the transpose of a matrix $\cM$,
$\diag(m_1,\dots,m_k)\in\bbC^{k\times k}$ a diagonal $k\times k$
matrix, and $\AC_{\loc}(\bbR)$ denotes the set of locally absolutely
continuous functions on $\bbR$. The spectrum, point spectrum (the
set of eigenvalues), essential spectrum, absolutely continuous spectrum,
and singularly continuous spectrum of a self-adjoint linear operator
$T$ in a separable complex Hilbert space are denoted by $\spec(T)$,
$\spec_{\p}(T)$, $\spec_{\ess}(T)$, $\spec_{\ac}(T)$, $\spec_{\singc}(T)$,
respectively.

The basic assumption for this section will be the following.

\begin{hypothesis}\lb{h2.1} ${}$ \\
$(i)$ Fix $m\in\bbN$, suppose $\cQ=\cQ^*\in L_{\loc}^1(\bbR)^{m\times
m}$ and introduce the differential expression
\begin{equation}
\cL=-\cI_m\f{d^2}{dx^2}+\cQ, \quad x\in\bbR. \lb{2.1}
\end{equation} 
$(ii)$ Suppose $\cL$ is in the limit point case at $\pm\infty$. 
\end{hypothesis}

Given Hypothesis~\ref{h2.1}\,(i) we consider the matrix-valued
Schr\"odinger equation 
\begin{equation}
-\psi^{\prime\prime}(z,x)+\cQ(x)\psi(z,x) =z\psi(z,x)  
\text{ for a.e. $x\in\bbR$},  \lb{2.6}
\end{equation}
where $z\in\bbC$ plays the role of a spectral parameter and $\psi$ is
assumed to satisfy
\begin{equation}
\psi(z,\cdot), \psi^\prime (z,\cdot) \in \AC_{\loc}(\bbR)^{m\times m}.
\lb{2.7}
\end{equation}
Throughout this paper, $x$-derivatives are abbreviated by a prime
$\prime$.

Let $\Psi(z,x,x_0)$ be a $2m\times 2m$ normalized fundamental system of 
solutions of \eqref{2.6} at some $x_0\in\bbR$ which we partition as
\begin{equation}
\Psi(z,x,x_0)
=\begin{pmatrix}\theta(z,x,x_0) & \phi(z,x,x_0)\\
\theta'(z,x,x_0)& \phi'(z,x,x_0)\end{pmatrix}. \lb{2.16}
\end{equation}
Here $\prime$ denotes $d/dx$, $\theta(z,x,x_0)$ and $\phi(z,x,x_0)$ are
$m\times m$ matrices, entire with respect to $z\in\bbC$, and normalized
according to 
\begin{equation}
\Psi(z,x_0,x_0)=\cI_{2m}, \lb{2.15}
\end{equation}
that is,
\begin{equation}
\theta(z,x_0,x_0)=\phi'(z,x_0,x_0)=\cI_m, \quad
\phi(z,x_0,x_0)=\theta'(z,x_0,x_0)=0. \lb{2.17}
\end{equation}
In this context, we briefly recall a set of formulas needed later in
Section \ref{s4}. Introducing 
\begin{equation} 
\cJ=\begin{pmatrix}0& -\cI_m \\ \cI_m & 0  \end{pmatrix},
\lb{2.15A}  
\end{equation}
one infers   
\begin{equation}\lb{2.310}
\Psi(\bar{z},x,x_0)^*\cJ\Psi(z,x,x_0)=\cJ,
\end{equation}
which implies $\cJ\Psi(z,x,x_0)(\Psi(\bar{z},x,x_0)\cJ)^*=\cI_{2m}$ 
and hence
\begin{equation}\lb{2.330}
\Psi(z,x,x_0)\cJ\Psi(\bar{z},x,x_0)^*=\cJ.
\end{equation}  
Writing out \eqref{2.310} and \eqref{2.330} explicitly yields 
\begin{align}
\theta'(\bar z,x,x_0)^*\theta(z,x,x_0)-
\theta(\bar z,x,x_0)^*\theta'(z,x,x_0)&=0, \lb{2.72} \\
\phi'(\bar z,x,x_0)^*\phi(z,x,x_0)-
\phi(\bar z,x,x_0)^*\phi'(z,x,x_0)&=0, \lb{2.73} \\
\phi'(\bar z,x,x_0)^*\theta(z,x,x_0)-
\phi(\bar z,x,x_0)^*\theta'(z,x,x_0)&=\cI_m, \lb{2.74} \\
\theta(\bar z,x,x_0)^*\phi'(z,x,x_0)-
\theta'(\bar z,x,x_0)^*\phi(z,x,x_0)&=\cI_m, \lb{2.75}
\end{align}
and 
\begin{align}
\phi(z,x,x_0)\theta(\bar z,x,x_0)^*-
\theta(z,x,x_0)\phi(\bar z,x,x_0)^*&=0, \lb{2.92} \\
\phi'(z,x,x_0)\theta'(\bar z,x,x_0)^*-
\theta'(z,x,x_0)\phi'(\bar z,x,x_0)^*&=0, \lb{2.93} \\
\phi'(z,x,x_0)\theta(\bar z,x,x_0)^*-
\theta'(z,x,x_0)\phi(\bar z,x,x_0)^*&=\cI_m, \lb{2.94} \\
\theta(z,x,x_0,\alpha)\phi'(\bar z,x,x_0)^*-
\phi(z,x,x_0,\alpha)\theta'(\bar z,x,x_0)^*&=\cI_m. \lb{2.95}
\end{align}
Next, assuming $ -\infty\le a< b \le \infty$, we consider the spaces 
\begin{equation}
N(z,\pm\infty)=\{\phi\in L^2((c,\pm\infty))^m \mid -\phi''+\cQ\phi=z\phi
\text{ a.e. on $(c,\pm\infty)$} \}, \lb{2.20A}
\end{equation}
for some $c\in\bbR$ and $z\in\bbC$. (Here
$(\phi,\psi)_{\bbC^n}=\sum_{j=1}^n \overline\phi_j\psi_j$
denotes the standard scalar product in $\bbC^n$, abbreviating
$\chi\in\bbC^n$ by $\chi=(\chi_1,\cdots,\chi_n)^t$.)  Both dimensions of
the spaces in \eqref{2.20A}, $\dim_\bbC(N(z,\infty))$ and
$\dim_\bbC(N(z,-\infty))$, are constant for $z\in\bbC_\pm=\{\zeta\in\bbC
\mid \pm\Im(\zeta)> 0 \}$ (see, e.g., \cite{KR74}). One then
recalls that $\cL$ in \eqref{2.1} is in the limit point case at
$\pm\infty$ whenever
\begin{equation}
\dim_\bbC(N(z,\pm\infty))=m \text{  for all
$z\in\bbC\backslash\bbR$.} \lb{2.20f}
\end{equation}
Since the potential $\cQ_\Sigma$ to be constructed in Section \ref{s4}
will automatically lead to the limit point case at $\pm\infty$, we decided
to limit  our considerations mainly to this situation. In this context we
note the well-known fact that if $\cL$ in \eqref{2.1} is in the limit
point case at $\pm\infty$, then the $m\times m$ Weyl--Titchmarsh matrices
associated with $\cL$, the half-lines $[x,\pm\infty)$, and a Dirichlet
boundary condition at $x$, are given by
\begin{equation}
\cM_{\pm}(z,x_0)=\Psi_{\pm}'(z,x,x_0)\Psi_\pm(z,x,x_0)^{-1}\big|_{x=x_0}, 
\quad z\in\bbC\backslash\bbR, \lb{2.47}
\end{equation}
where $\Psi_{\pm}$ satisfy
$(\cL-z\cI_m)\Psi_{\pm}(z,\cdot,x_0)=0$ and 
\begin{equation}
\Psi_\pm(z,\cdot,x_0)\in L^2([x_0,\pm\infty))^{m\times m}. \lb{2.20g}
\end{equation}
The actual normalization of $\Psi_\pm(z,\cdot,x_0)$ is clearly irrelevant
and hence $\Psi_\pm(z,\cdot,x_0)$ can be replaced by $\Psi_\pm
(z,\cdot,x_0)C$, where $C$ is any nonsingular $m\times m$ matrix. 

For later reference we summarize the principal results on
$\cM_{\pm}(z,x_0)$ in the next theorem. First we recall the 
following definition.

\begin{definition} \lb{d2.6}
A map $\cM\colon\bbC_+\to \bbC^{n\times n}$, $n\in\bbN$, extended to
$\bbC_-$ by $\cM(\bar z)= \cM(z)^*$ for all $z\in\bbC_+$, is called an
$n\times n$ Herglotz matrix\footnote{There appears to be 
considerable confusion in the literature since Nevanlinna, Pick,
Nevanlinna--Pick matrix, in addition to 
Herglotz matrix, are also in use.  In part these discrepancies can be
traced back to  the use of the upper half-plane $\bbC_+$ versus the open
unit disk $D$, and in some cases the geographical location of the author
in question determines the preferred notation. Following a 
tradition in mathematical physics, we adopt the notion of Herglotz
functions in this paper.} if it is analytic on $\bbC_+$ and
$\Im(\cM(z))\ge 0$ for all $z\in\bbC_+$.
\end{definition}

\noindent Here we denote $\Im(\cM)=(\cM-\cM^*)/2i$ and
$\Re(\cM)=(\cM+\cM^*)/2$.

In the scalar context $n=1$, the condition $\Im(\cM(z))\geq 0$ in
Definition \ref{d2.6} can be replaced by $\Im(m(z)) > 0$ for the
corresponding scalar counterpart $m(z)$.

\begin{theorem}
[\cite{CG01}, \cite{GT00}, \cite{HS81}, 
\cite{HS82}, \cite{HS86}, \cite{KS88}] \lb{t2.7} 
Assume Hypothesis~\ref{h2.1} and suppose that 
$z\in\bbC\backslash\bbR$, $x_0\in\bbR$. Then \\
$(i)$ $\pm \cM_{\pm}(z,x_0)$ is an $m\times m$ matrix-valued 
Herglotz function of maximal rank. In particular,
\begin{gather}
\Im(\pm \cM_{\pm}(z,x_0)) > 0, \quad z\in\bbC_+, 
\lb{2.21} \\
\cM_{\pm}(\overline z,x_0)=\cM_{\pm}( z,x_0)^*,\lb{2.22} \\
\rank (\cM_{\pm}(z,x_0))=m, \lb{2.23} \\
\lim_{\varepsilon\downarrow 0} \cM_{\pm}(\lambda+
i\varepsilon,x_0) \text{
exists for a.e.\
$\lambda\in\bbR$}. \lb{2.24a}
\end{gather}
Local singularities of $\pm \cM_{\pm}(z,x_0)$ and 
$\mp \cM_{\pm}(z,x_0)^{-1}$ are necessarily real and at most of first
order in the sense that 
\begin{align}
&\mp \lim_{\epsilon\downarrow0}
\left(i\epsilon\,
\cM_{\pm}(\lambda+i\epsilon,x_0)\right) \geq 0, \quad \lambda\in\bbR,
\lb{2.24b} \\ 
& \pm \lim_{\epsilon\downarrow0}
i\epsilon\cM_{\pm}(\lambda+i\epsilon,x_0)^{-1}
\geq 0, \quad \lambda\in\bbR. \lb{2.24c}
\end{align}
$(ii)$  $\pm \cM_{\pm}(z,x_0)$ admit the representations
\begin{equation}
\pm \cM_{\pm}(z,x_0)=\Re(\pm \cM_{\pm}(\pm i,x_0)) +\int_\bbR 
d\Omega_\pm(\lambda,x_0) \,
\big((\lambda-z)^{-1}-\lambda(1+\lambda^2)^{-1}\big), \lb{2.25a} 
\end{equation}
where
\begin{equation}
\int_\bbR \|d\Omega_\pm(\lambda,x_0)\|_{\bbC^{m\times m}}\,
(1+\lambda^2)^{-1}<\infty \lb{2.27} 
\end{equation}
and
\begin{equation}
\Omega_\pm((\lambda,\mu],x_0)=\lim_{\delta\downarrow
0}\lim_{\varepsilon\downarrow 0}\f1\pi
\int_{\lambda+\delta}^{\mu+\delta} d\nu \, \Im(\pm
\cM_\pm(\nu+i\varepsilon,x_0)). \lb{2.29} 
\end{equation}
$(iii)$  Define the $2m\times m$ matrices
\begin{align}
\Psi_\pm(z,x,x_0)&=\begin{pmatrix}\psi_{\pm}(z,x,x_0)\\
\psi_{\pm}' (z,x,x_0)  \end{pmatrix} \no \\
&=\begin{pmatrix}\theta(z,x,x_0) & \phi(z,x,x_0)\\
\theta'(z,x,x_0)& \phi'(z,x,x_0)\end{pmatrix} 
\begin{pmatrix} \cI_m \\
\cM_\pm(z,x_0) \end{pmatrix}, \lb{2.31}
\end{align}
then 
\begin{equation}
\Im(\cM_\pm(z,x_0))=\Im(z) \int_{x_0}^{\pm\infty}dx \,
\psi_\pm(z,x,x_0)^* \psi_\pm(z,x,x_0). \lb{2.32}
\end{equation}
$(iv)$ Denote by $C_\varepsilon\subset\bbC_+$ the open sector with vertex
at zero, symmetry axis along the positive imaginary axis, and opening
angle $\varepsilon$, with $0<\varepsilon<\pi/2$. Then
\begin{equation}
\cM_\pm(z,x_0)\underset{\substack{|z|\to\infty\\z\in C_\varepsilon}}{=}
\pm z^{1/2} \cI_m+ o(1). \lb{2.32a}
\end{equation}
\end{theorem}

Necessary and sufficient conditions for $\cM_\pm(z,x_0)$ to be the
half-line $m\times m$ Weyl--Titchmarsh matrix associated with a
Schr\"odinger operator on $[x_0,\pm\infty)$ in terms of the
corresponding measures $\Omega_\pm(\cdot,x_0)$ in the Herglotz
representation \eqref{2.25a} of $\cM_\pm(z,x_0)$ can be derived using
the matrix-valued extension of the classical inverse spectral theory
approach due to Gelfand and Levitan \cite{GL51}, as worked out by
Rofe-Beketov \cite{Ro60}. The following result describes sufficient
conditions for a monotonically nondecreasing matrix function to be the
matrix spectral function of a half-line Schr\"odinger operator. It extends
well-known results in the scalar case $m=1$ (cf.\
\cite[Sects.\ 2.5, 2.9]{Le87}, \cite{LG64}, \cite[Sect.\ 26.5]{Na68},
\cite{Th79}).

\begin{theorem} [\cite{Ro60}] \lb{t2.7a} 
Let $\Omega_+$ be a monotonically nondecreasing $m\times m$
matrix-valued function on $\bbR$ satisfying the following two
conditions. \\
\noindent $(i)$ 
Whenever $f\in C([x_0,\infty))^{m\times 1}$ with compact support 
contained in $[x_0,\infty)$ and
\begin{equation}
\int_\bbR F(\lambda)^*d\Omega_+(\lambda)\,F(\lambda) =0, 
 \text{ then $f=0$~a.e.,} \lb{2.35}
\end{equation} 
where 
\begin{equation}
F(\lambda)=\lim_{R\uparrow\infty}\int_{x_0}^R
dx\,\f{\sin(\lambda^{1/2}(x-x_0))}{\lambda^{1/2}}f(x), \quad
\lambda\in\bbR.
\lb{2.36}
\end{equation}
\noindent $(ii)$ Define 
\begin{equation}
\wti\Omega_+(\lambda)=\begin{cases}
\Omega_+(\lambda)-\f{2}{3\pi}\lambda^{3/2}, & \lambda\geq 0 \\
\Omega_+(\lambda), & \lambda<0 \end{cases} \lb{2.37}
\end{equation}
and assume the limit 
\begin{equation}
\lim_{R\uparrow\infty}\int_{-\infty}^R d\wti\Omega_+(\lambda) \,
\f{\sin(\lambda^{1/2}(x-x_0))}{\lambda^{1/2}}= \Phi(x) \lb{2.38}
\end{equation}
exists and $\Phi\in L^\infty([x_0,R])^{m\times m}$ for all $R>x_0$.
Moreover, suppose that for some $r\in\bbN_0$, $\Phi^{(r+1)}\in
L^1([x_0,R])^{m\times m}$ for all $R>x_0$, and $\Phi(x_0)=0$. \\
Then $\Omega_+$ is the matrix spectral function of a self-adjoint
Schr\"odinger operator $H_+$ in $L^2([x_0,\infty))^{m}$ associated
with the $m\times m$ matrix-valued differential expression 
$\cL_+=-d^2/dx^2 \cI_m+\cQ$, $x>x_0$, with a Dirichlet boundary condition
at $x_0$, a self-adjoint boundary condition at $\infty$ $($if
necessary$)$, and a self-adjoint potential matrix $\cQ$ with
$\cQ^{(r)}\in L^1([x_0,R])^{m\times m}$ for all $R>x_0$. 
\end{theorem}

Next, assuming Hypothesis \ref{h2.1}, we introduce the self-adjoint 
Schr\"odinger operator $H$ in $L^2(\bbR)^m$ by
\begin{align}
&H=-\cI_m \f{d^2}{dx^2}+\cQ, \lb{2.39} \\
&\dom(H)=\{g\in L^2(\bbR)^m \mid g,g^\prime\in
\AC_{\loc}(\bbR)^m;\, (-g^{\prime\prime}+\cQ g)\in L^2(\bbR)^m\}. \no
\end{align}
The resolvent of $H$ then reads 
\begin{equation}
((H-z)^{-1}f)(x)= \int_\bbR dx^\prime\, \cG(z,x,x^\prime)f(x^\prime), 
\quad z\in\bbC\backslash\bbR, \; f\in L^2(\bbR)^m, \lb{2.38a} 
\end{equation}
with the Green's matrix $\cG(z,x,x')$ of $H$ given by 
\begin{align}
\cG(z,x,x^\prime)=\psi_\mp(z,x,x_0)[\cM_-(z,x_0) &
-\cM_+(z,x_0)]^{-1}\psi_\pm(\overline z,x^\prime,x_0)^*, 
\no \\
& \hspace*{2.05cm} x\lesseqgtr x^\prime,\; 
z\in\bbC\backslash\bbR. \lb{2.33}
\end{align}
Introducing
\begin{equation}
\cN_\pm (z,x_0)=\cM_-(z,x_0)\pm \cM_+(z,x_0),
\quad z\in\bbC\backslash\bbR, \lb{2.33A} 
\end{equation}
the $2m\times 2m$ Weyl--Titchmarsh function $\cM(z,x_0)$ associated with
$H$ on $\bbR$ is then given by 
\begin{align}
\cM(z,x_0)&=\big(\cM_{j,j^\prime}(z,x_0)\big)_{j,j^\prime=1,2} \lb{2.34}
\\ &=\begin{pmatrix}
\cM_\pm (z,x_0)\cN_-(z,x_0)^{-1}\cM_\mp (z,x_0)   
&\cN_-(z,x_0)^{-1}\cN_+(z,x_0)/2  \\ 
\cN_+(z,x_0)\cN_-(z,x_0)^{-1}/2  &
\cN_-(z,x_0)^{-1} \end{pmatrix},  \no \\
& \hspace*{8.95cm}  z\in\bbC\backslash\bbR. \no
\end{align} 

The basic results on $\cM(z,x_0)$ then read as follows.

\begin{theorem} [\cite{GT00}, \cite{HS81}, \cite{HS82}, 
\cite{HS86}, \cite{KS88}] \lb{t2.8} 
Assume Hypothesis~\ref{h2.1} and suppose that $z\in\bbC 
\backslash \bbR$, $x_0\in\bbR$.  Then, \\
$(i)$ $\cM(z,x_0)$ is a matrix-valued Herglotz function of rank 
$2m$ with
representation
\begin{equation}
\cM(z,x_0)=\Re(\cM(i,x_0)) +\int_\bbR d\Omega(\lambda,x_0)\,
\big((\lambda-z)^{-1}-\lambda(1+\lambda^2)^{-1}\big), \lb{2.42} 
\end{equation}
where
\begin{equation}
\int_\bbR \Vert d\Omega(\lambda,x_0)
\Vert_{\bbC^{2m\times 2m}} \,(1+\lambda^2)^{-1}<\infty \lb{2.44}
\end{equation}
and
\begin{equation}
\Omega((\lambda,\mu],x_0)=\lim_{\delta\downarrow
0}\lim_{\varepsilon\downarrow 0}\f1\pi
\int_{\lambda+\delta}^{\mu+\delta} d\nu \, 
\Im(\cM(\nu+i\varepsilon,x_0)). \lb{2.46}
\end{equation}
$(ii)$ $z\in\bbC\backslash\spec(H)$ if and only if $\cM(z,x_0)$ is
holomorphic  near $z$. 
\end{theorem}
Here $\spec(T)$ denotes the spectrum of $T$. Later on we will denote by 
$\spec_{\text{ac}}(T)$ the absolutely continuous spectrum of $T$.

Finally, we state the following characterization of $\cM(z,x_0)$ to be 
used later on. In the scalar context $m=1$ this has been used by
Rofe-Beketov \cite{Ro67}, \cite{Ro91} (see also \cite[Sect.\ 7.3]{Le87}).

\begin{theorem} [\cite{Ro67}, \cite{Ro91}] \lb{t2.9} 
Assume Hypothesis~\ref{h2.1}, suppose that $z\in\bbC 
\backslash \bbR$, $x_0\in\bbR$, and let $\ell, r\in\bbN_0$. Then the
following assertions are equivalent. \\
$(i)$ $\cM(z,x_0)$ is the $2m\times 2m$ 
Weyl--Titchmarsh matrix associated with a Schr\"odinger operator 
$H$ in $L^2(\bbR)^m$ of the type \eqref{2.39} with an $m\times m$
matrix-valued potential $\cQ\in L^1_{\loc}(\bbR)$ and $\cQ\in
C^\ell((-\infty,x_0))$ and $\cQ\in C^r((x_0,\infty))$. \\
$(ii)$ $\cM(z,x_0)$ is of the type \eqref{2.34} with $\cM_\pm(z,x_0)$
being half-line $m\times m$ Weyl--Titchmarsh matrices on $[x_0,\pm\infty)$
corresponding to a Dirichlet boundary condition at $x_0$ and a
self-adjoint boundary condition at $-\infty$ and/or $\infty$ $($if
any$)$ which are associated with an $m\times m$ matrix-valued potential
$\cQ$ satisfying $\cQ\in C^\ell((-\infty,x_0))$ and $\cQ\in
C^r((x_0,\infty))$, respectively. \\ 
If $(i)$ or $(ii)$ holds, then the $2m\times 2m$ matrix-valued spectral
measure $\Omega(\cdot,x_0)$ associated with $\cM(z,x_0)$ is determined
by \eqref{2.34} and \eqref{2.46}.
\end{theorem}

Next, we consider variations of the reference point $x_0\in\bbR$. In
analogy to \eqref{2.47}, we note that in the case where the Schr\"odinger
differential expression $\cL$ is in the limit point case at $\pm\infty$, 
\begin{equation}
\cM_{\pm}(z,x)=\Psi_{\pm,x}'(z,x,x_0)\Psi_\pm(z,x,x_0)^{-1}, 
\quad z\in\bbC\backslash\bbR, \lb{2.48}
\end{equation}
represents the corresponding half-line Weyl--Titchmarsh matrix on
$[x,\pm\infty)$, $x\in\bbR$, with $\Psi_\pm(z,\cdot,x_0)$ defined in
\eqref{2.31}. Again the actual normalization of $\Psi_\pm$ is, of course,
irrelevant. Since $\Psi_\pm$ satisfies the second-order linear $m\times
m$ matrix-valued differential equation
\eqref{2.6}, $\cM_\pm$ in \eqref{2.48} satisfies the matrix-valued
Riccati-type equation (independently of any limit point assumptions at
$\pm\infty$)
\begin{equation}
\cM_{\pm}'(z,x)+\cM_\pm(z,x)^2=\cQ(x)-z \cI_m, \quad x\in\bbR,
\;  z\in\bbC\backslash\bbR. \lb{2.49}
\end{equation}

The asymptotic high-energy behavior of $\cM_\pm(z,x)$ as $|z|\to\infty$
has recently been determined in \cite{CG01} under minimal smoothness
conditions on $\cQ$ and without assuming that $\cL$ is in the limit point
case at $\pm\infty$. Here we recall just a special case of the asymptotic
expansion proved in \cite{CG01} which is most suited for our discussion in
Section \ref{s4}. We denote by $C_\varepsilon\subset\bbC_+$ the open
sector with vertex at zero, symmetry axis along the positive imaginary
axis, and opening angle $\varepsilon$, with $0<\varepsilon< \pi/2$.

\begin{theorem} [\cite{CG01}] \lb{t2.10}
Fix $x_0\in\bbR$ and let $x\geq x_0$. In addition to
Hypothesis~\ref{h2.1} suppose that 
$\cQ\in C^\infty([x_0,\pm\infty))^{m\times m}$ and that $\cL$ is in the
limit point case at $\pm\infty$. Let $\cM_\pm(z,x)$, $x\geq x_0$, be
defined as in \eqref{2.48}. Then, as $\abs{z}\to\infty$ in
$C_\varepsilon$, $\cM_\pm(z,x)$ has an asymptotic expansion of the form
$(\Im(z^{1/2})>0$, $z\in\bbC_+)$
\begin{equation}
\cM_\pm(z,x)\underset{\substack{\abs{z}\to\infty\\ z\in
C_\varepsilon}}{=} \pm i \cI_m z^{1/2}+\sum_{k=1}^N
\cM_{\pm,k}(x)z^{-k/2}+ o(|z|^{-N/2}), \quad N\in\bbN. \lb{2.50}
\end{equation}
The expansion \eqref{2.50} is uniform with respect to $\arg\,(z)$ for 
$|z|\to \infty$ in $C_\varepsilon$ and uniform in $x$ as long as $x$
varies in compact subsets of $[x_0,\infty)$. The expansion coefficients
$\cM_{\pm,k}(x)$ can be recursively computed from 
\begin{align}
\begin{split}
\cM_{\pm,1}(x)&=\mp\f{i}{2} \cQ(x),
\quad \cM_{\pm,2}(x)= \f1{4} \cQ^\prime(x),  \\
\cM_{\pm,k+1}(x)&=\pm\f{i}2\bigg(\cM_{\pm,k}^\prime(x)+
\sum_{\ell=1}^{k-1}\cM_{\pm,\ell}(x) \cM_{\pm,k-\ell}(x) \bigg),
\quad k\ge 2. \lb{2.51}
\end{split}
\end{align} 
The asymptotic expansion \eqref{2.50} can be differentiated to any
order with respect to $x$. 
\end{theorem}

\begin{remark} \lb{r2.11} ${}$ \\
$(i)$ Due to the recursion relation \eqref{2.51}, the coefficients
$\cM_{\pm,k}$ are universal polynomials in $\cQ$ and its $x$-derivatives
$($i.e., differential polynomials in $\cQ$$)$. That the asymptotic
expansion
\eqref{2.50} can be differentiated to arbitrary order in $x$ follows from
repeated use of the Riccati-type equation \eqref{2.49}. \\
$(ii)$ In the case where $\cQ$ and its $x$-derivatives are in
$L^1(\bbR)^{m\times m}$, or in the case where $\cQ$ is periodic and hence
Floquet theory applies, the proof of the existence of an asymptotic
expansion of the type \eqref{2.50} follows in a routine manner by
iterating appropriate Volterra-type integral equations. The general case,
however, is intricate as is evident from the treatment in \cite{CG01}. 
\end{remark}

Finally, in addition to \eqref{2.33}, one infers for the $2m\times 2m$ 
Weyl--Titchmarsh function $\cM(z,x)$ associated with $H$ on $\bbR$ in
connection with arbitrary half-lines $[x,\pm\infty)$, $x\in\bbR$, 
\begin{align}
\cM(z,x)&=\big(\cM_{j,j^\prime}(z,x)\big)_{j,j^\prime=1,2},
\quad  z\in\bbC\backslash\bbR,  \lb{2.52} \\
\cM_{1,1}(z,x)&=\cM_\pm(z,x)[\cM_-(z,x)-\cM_+(z,x)]^{-1}
\cM_\mp(z,x) \no \\ 
&=\psi_+'(z,x,x_0)[\cM_-(z,x_0)-\cM_+(z,x_0)]^{-1}\psi_-'(\ol z,x,x_0)^*,
\lb{2.53} \\
\cM_{1,2}(z,x)&=2^{-1} [\cM_-(z,x)-\cM_+(z,x)]^{-1}
[\cM_-(z,x)+\cM_+(z,x)] \no \\
&=\psi_+(z,x,x_0)[\cM_-(z,x_0)-\cM_+(z,x_0)]^{-1}\psi_-'(\ol z,x,x_0)^*,
\lb{2.54} \\
\cM_{2,1}(z,x)&=2^{-1}
[\cM_-(z,x)+\cM_+(z,x)][\cM_-(z,x)-\cM_+(z,x)]^{-1} \no \\
&=\psi_+'(z,x,x_0)[\cM_-(z,x_0)-\cM_+(z,x_0)]^{-1}\psi_-(\ol z,x,x_0)^*,
\lb{2.55} \\
\cM_{2,2}(z,x)&=[\cM_-(z,x)-\cM_+(z,x)]^{-1} \no \\
&=\psi_+(z,x,x_0)[\cM_-(z,x_0)-\cM_+(z,x_0)]^{-1}\psi_-(\ol z,x,x_0)^*.
\lb{2.56} 
\end{align}
Introducing the convenient abbreviation,
\begin{equation}
\cM(z,x)=\begin{pmatrix} \gh(z,x) & -\gg_2(z,x) \\ -\gg_1(z,x) & \gg(z,x)
\end{pmatrix}, \quad z\in\bbC\backslash\bbR, \; x\in\bbR, \lb{2.57}
\end{equation}
one then verifies from \eqref{2.52}--\eqref{2.57} and from $\cM(\ol
z,x)^*=\cM(z,x)$, $\cM_\pm(\ol z,x)^*=\cM_\pm(z,x)$ that 
\begin{align}
&\gg(\ol z,x)^*=\gg(z,x), \quad \gg_2(\ol z,x)^*=\gg_1(z,x), 
\quad \gh(\ol z,x)^*=\gh(z,x), \lb{2.58} \\
&\gg(z,x)\gg_1(z,x)=\gg_2(z,x)\gg(z,x), \lb{2.59} \\
&\gh(z,x)\gg_2(z,x)=\gg_1(z,x)\gh(z,x), \lb{2.59a} \\
&\gg(z,x)=[\cM_-(z,x)-\cM_+(z,x)]^{-1}, \lb{2.60} \\
&\gg(z,x)\gh(z,x)-\gg_2(z,x)^2=-(1/4)\cI_m, \lb{2.61} \\
&\gh(z,x)\gg(z,x)-\gg_1(z,x)^2=-(1/4)\cI_m, \lb{2.62} \\ 
&\cM_\pm(z,x)=\mp(1/2)\gg(z,x)^{-1}-\gg(z,x)^{-1}\gg_2(z,x)  \lb{2.63} \\
&\qquad \qquad \;\;\;\,\, \,\mp(1/2)\gg(z,x)^{-1}-\gg_1(z,x)\gg(z,x)^{-1}.
\lb{2.64}
\end{align}
Moreover, the Riccati-type equations \eqref{2.49} imply the following
results needed in Section \ref{s4}.

\begin{lemma} \lb{l2.12}
Let $z\in\bbC\backslash\bbR$ and define $\cM_\pm$ by
\eqref{2.48} so that $\cM_\pm$ satisfy the Riccati-type equation
\eqref{2.49}. Then, for a.e.\ $x\in\bbR$,
\begin{align}
\gg'&=-(\gg_1+\gg_2), \lb{2.65} \\
\gg_1'&=-(\cQ-z\cI_m)\gg-\gh \lb{2.66} \\
& = (-\gg''+\gg\cQ-\cQ\gg)/2, \lb{2.66a} \\ 
\gg_2'&=-\gg(\cQ-z\cI_m)-\gh \lb{2.67} \\
& =(-\gg''+\cQ\gg-\gg\cQ)/2, \lb{2.67a} \\
\gh'&=-\gg_1(\cQ-z\cI_m)-(\cQ-z\cI_m)\gg_2, \lb{2.68} \\  
\gh&=[\gg''-\gg(\cQ-z\cI_m)-(\cQ-z\cI_m)\gg]/2 \lb{2.69} 
\end{align}
if $\cQ\in L^1_{\loc}(\bbR)^{m\times m}$, and 
\begin{align}
\gg_1''&=-2(\cQ-z\cI_m)\gg'-\cQ'\gg+\gg_1\cQ-\cQ\gg_1, \lb{2.66b} \\
\gg_2''&=-2\gg'(\cQ-z\cI_m)-\gg\cQ'+\cQ\gg_2-\gg_2\cQ \lb{2.67b}  
\end{align}
if in addition $\cQ'\in L^1_{\loc}(\bbR)^{m\times m}$.
\end{lemma}
\begin{proof}
\eqref{2.49} rewritten in terms of $\gg$, $\gg_1$, $\gg_2$ yields
\begin{align}
&\pm(1/2)\gg^{-1}\gg'\gg^{-1}+\gg^{-1}\gg'\gg^{-1}\gg_2-\gg^{-1}\gg_2'
+(1/4)\gg^{-2}+\gg^{-1}\gg_2\gg^{-1}\gg_2 \no \\
&\pm(1/2)\gg^{-2}\gg_2 \pm(1/2)\gg^{-1}\gg_2\gg^{-1}=\cQ-z\cI_m. \lb{2.70}
\end{align}
Taking the difference of the two equations in \eqref{2.70} yields
\eqref{2.65}. Adding the two equations in \eqref{2.70} and using
\eqref{2.59}, \eqref{2.61}, \eqref{2.62}, and \eqref{2.65} yields
\eqref{2.66} and \eqref{2.67}. Combining \eqref{2.66}, \eqref{2.67}, and
\eqref{2.65} implies \eqref{2.69}. Inserting \eqref{2.69} into
\eqref{2.66} and \eqref{2.67} yields \eqref{2.66a} and \eqref{2.67a}.
\eqref{2.68} follows from differentiating
$\gh=\gg_1^2\gg^{-1}-(1/4)\gg^{-1}$, inserting
$\gg'$ and $\gg_1'$ from \eqref{2.65} and \eqref{2.66}, and making
repeated use of the identities \eqref{2.59}, \eqref{2.62}. Finally,
\eqref{2.66b} (resp.\ \eqref{2.67b}) follows from differentiating
\eqref{2.66} (resp.\ \eqref{2.67}) inserting \eqref{2.68} for $\gh'$.
Alternatively, \eqref{2.67}--\eqref{2.67b} follow directly from
\eqref{2.66}--\eqref{2.66b} using \eqref{2.58}.
\end{proof}

\section{Polynomial Pencils of Matrices in a Nutshell}
\lb{s3}

Since self-adjoint polynomial pencils of matrices play a role in our
principal section \ref{s4}, we briefly review some of the corresponding
definitions and basic results, mainly following the monograph of Markus
\cite{Ma88a} and papers by Markus and Matsaev \cite{MM76},
\cite{MM85}. While all results below are discussed for operator pencils
by Markus and Matsaev, we will only quote them in the matrix context,
for simplicity.

Given $m\in\bbN$, we denote by
\begin{equation}
\cA(z)=\sum_{k=0}^n \cA_k z^k, \quad \cA_k\in\bbC^{m\times m}, 
\, 1\leq k\leq n, \; z\in\bbC, \lb{3.1}
\end{equation}
a polynomial pencil of $m\times m$ matrices (in short, a pencil) in the
following. $\cA$ is called  of degree $n\in\bbN_0$ if $\cA_n\neq 0$
and monic if $\cA_n=\cI_m$.

\begin{definition} \lb{d3.1} Let $\cA$ be a pencil of the type 
\eqref{3.1}. \\
$(i)$ The {\it spectrum} of $\cA$, denoted by $\spec(\cA)$, is defined
by 
\begin{equation}
\spec(\cA)=\{z\in\bbC\,|\, \cA(z) \;\text{is not invertible}\}. 
\end{equation}
$z_0\in\bbC$ is called an eigenvalue of $\cA$ if $\cA(z_0)f_0=0$ has a
solution $f_0\in\bbC^m\backslash\{0\}$. \\   
$(ii)$ A monic pencil $\cC$ is called a {\it $($right\,$)$ divisor}
of $\cA$ if 
\begin{equation}
\cA(z)=\cB(z)\cC(z), \quad z\in\bbC
\end{equation}
for some pencil $\cB$. If in addition
$\spec(\cB)\cap\spec(\cC)=\emptyset$, then
$\cC$ is called a {\it $($right\,$)$ spectral divisor} of $\cA$. \\
$(iii)$ $\cZ\in\bbC^{m\times m}$ is called a {\it $($matrix\,$)$ root}
of the pencil $\cA$ if 
\begin{equation}
\cA(\cZ)=0,
\end{equation}
where $\cA(\cZ)$ is defined as
\begin{equation}
\cA(\cZ)=\sum_{k=0}^n \cA_k \cZ^k.
\end{equation}
$\cZ\in\bbC^{m\times m}$ is called a {\it $($matrix\,$)$ spectral root}
of the pencil $\cA$ if $(z\cI_m-\cZ)$ is a spectral divisor of
$\cA$. \\
$(iv)$ The pencil $\cA$ is called {\it self-adjoint} if
$\cA_k=\cA_k^*$ for all $1\leq k\leq n$ $($i.e., $\cA(\ol z)^*=\cA(z)$ for
all $z\in\bbC$$)$. \\
$(v)$ A self-adjoint pencil $\cA$ is called {\it weakly hyperbolic} if
$\cA_n>0$ and for all $f\in\bbC^m\backslash\{0\}$, the roots of the
polynomial $(f,\cA(\cdot)f)_{\bbC^m}$ are real. If in addition all these
zeros are distinct, the pencil $\cA$ is called {\it hyperbolic}. \\
$(vi)$ Let $\cA$ be a weakly hyperbolic pencil and denote by 
$\{p_j(\cA,f)\}_{1\leq j\leq n}$,
\begin{equation}
 p_j(\cA,f)\leq
p_{j+1}(\cA,f), \; 1\leq j\leq n-1, \; f\in\bbC^m\backslash\{0\},
\end{equation}
the roots of the polynomial $(f,\cA(\cdot)f)_{\bbC^m}$ ordered in
magnitude. The range of the roots 
$p_j(\cA,f)$, $f\in\bbC^m\backslash\{0\}$ is denoted by
$\Delta_j(\cA)$ and called the {\it $j$th root zone} of $\cA$. \\
$(vii)$ A hyperbolic pencil $\cA$ is called {\it strongly hyperbolic}
if $\ol {\Delta_j(\cA)}$ and $\ol {\Delta_k(\cA)}$ are mutually
disjoint for $j\neq k$, $1\leq j,k \leq n$. 
\end{definition}

Moreover, the Vandermonde matrix corresponding to a collection
$\{\cZ_1,\dots,\cZ_n\}\subset\bbC^{m\times m}$ is defined by 
\begin{equation}
\cV(\cZ_1,\dots,\cZ_n)=\begin{pmatrix} \cI_m &
\cI_m & \dots & \cI_m \\ \cZ_1 & \cZ_2 & \dots & \cZ_n \\
\vdots & \vdots & \ddots & \vdots \\
\cZ_1^{n-1} & \cZ_2^{n-1} & \dots & \cZ_n^{n-1} \end{pmatrix}.
\end{equation}

\begin{theorem} [Markus \cite{Ma88a}, Sect.\ 29] \lb{t3.2} Let $\cA$ be
a pencil of the type \eqref{3.1}. \\
$(i)$ $(z\cI_m-\cZ)$ is a divisor of $\cA$ if and only if
$\cA(\cZ)=0$. $($This justifies the notation introduced in the last
part of  Definition \ref{d3.1}\,$(ii).)$\\
$(ii)$ Let $\cA$ be a monic pencil of degree $n$ and
$\cZ_1,\dots,\cZ_n$ spectral roots of $\cA$. Then the following
assertions are equivalent: \\
\indent $(\alpha)$ $\cV(\cZ_1,\dots,\cZ_n)$ is invertible. \\
\indent $(\beta)$ $\spec(\cZ_j)\cap\spec(\cZ_k)=\emptyset$, $j\neq k$,
$1\leq j,k\leq n$. \\
\indent $(\gamma)$ $\spec(\cA)=\bigcup_{j=1}^n \spec(\cZ_j)$. \\ 
$(iii)$ If $\cA$ is a self-adjoint pencil, then $\spec(\cA)$ is
symmetric with respect to $\bbR$. 
\end{theorem}

\begin{theorem} [Markus \cite{Ma88a}, Sect.\ 31] \lb{t3.3} Let $\cA$ be
a pencil of the type \eqref{3.1}. \\
$(i)$ If $\cA$ is a weakly hyperbolic pencil, then
$\spec(\cA)\subseteq\bbR$. \\
$(ii)$ The root zones $\Delta_j(\cA)$, $1\leq j\leq n$, of a weakly
hyperbolic pencil $\cA$ are intervals $($possibly degenerating to a
point\,$)$. \\
$(iii)$ If $\cA$ is a weakly hyperbolic pencil of order $n$ with root
zones $\{\Delta_j(\cA)\}_{1\leq j\leq n}$ and
$\Delta_k(\cA)=\{\lambda_k\}$ for some $k\in\{1,\dots,n\}$, then 
$\cA(z)=(z-\lambda_k)\cB(z)$, where $\cB$ is a weakly hyperbolic
pencil of order $n-1$, with root zones 
\begin{equation} 
\Delta_1(\cA),\dots,\Delta_{k-1}(\cA),\Delta_{k+1}(\cA),
\dots,\Delta_n(\cA). 
\end{equation}
$(iv)$ If $\cA$ is a weakly hyperbolic pencil then
$\Delta_j(\cA)\cap\Delta_{j+1}(\cA)$, $1\leq j\leq n-1$, consists of
at most one point. \\
$(v)$ If $\cA$ is a hyperbolic pencil, then 
$\Delta_j(\cA)\cap\Delta_k(\cA)=\emptyset$, $j\neq k$, $1\leq j,k\leq
n$.  Thus, a hyperbolic pencil is strongly hyperbolic if and only if 
\begin{equation}
\ol {\Delta_j(\cA)}\cap\ol {\Delta_{j+1}(\cA)}=\emptyset, 
\quad 1\leq j\leq n-1.
\end{equation}
$(vi)$ Suppose $\cA$ is a self-adjoint pencil of degree $n$,
$\cA_n>0$, and $\cA(\lambda)\neq 0$ for all $\lambda\in\bbR$. Then
$\cA$ is a weakly hyperbolic pencil if and only if there exist numbers
$\gamma_1<\gamma_2<\dots<\gamma_{n-1}$ such that $(-1)^j
\cA(\gamma_j)\geq 0$, $1\leq j\leq n-1$. \\
$(vii)$ Suppose $\cA$ is a self-adjoint pencil of degree $n$ with
$\cA_n>0$. Then $\cA$ is a strongly hyperbolic pencil if and only if
there exist numbers $\gamma_1<\gamma_2<\dots<\gamma_{n-1}$ such that
$(-1)^j \cA(\gamma_j) > 0$, $1\leq j\leq n-1$.
\end{theorem} 

\begin{theorem} [Markus \cite{Ma88a}, Sect.\ 31] \lb{t3.4} Let $\cA$ be
a pencil of the type \eqref{3.1}. \\
$(i)$ Suppose $\cA$ is a weakly hyperbolic pencil and 
$\ol {\Delta_{j_{0}-1}(\cA)}\cap\ol {\Delta_{j_{0}}(\cA)}=\ol
{\Delta_{j_{0}}(\cA)}\cap\ol {\Delta_{j_{0}+1}(\cA)}=\emptyset$ for
some $j_0 \in\{1,\dots,n\}$. Then $\cA$ has a spectral root
$\cZ_{j_0}$ such that $\spec(\cA)\cap\ol
{\Delta_{j_0}(\cA)}=\spec(\cZ_{j_0})$ and
$\cZ_{j_0}$ is similar to a self-adjoint matrix. \\
$(ii)$ A strongly hyperbolic pencil $\cA$ has $n$ spectral roots
$\{\cZ_j\}_{1\leq j\leq n}$ such that $\spec(\cZ_j)$
$=\spec(\cA)\cap\ol {\Delta_j(\cA)}$ and each $\cZ_j$, $1\leq j\leq
n$, is similar to a self-adjoint matrix.
\end{theorem} 

\begin{theorem} [Markus \cite{Ma88a}, Sect.\ 31, Markus and Matsaev
\cite{MM76}, \cite{MM85}] \lb{t3.5} Let $\cA$ be a pencil of the type
\eqref{3.1}. \\
$(i)$ A weakly hyperbolic monic pencil $\cA$ of degree $n$ is
decomposable as 
\begin{equation}
\cA(z)=(z\cI_m-\cY_n)(z\cI_m-\cY_{n-1})\cdots
(z\cI_m-\cY_1),
\end{equation}
with $\spec(\cY_j)\subset \ol {\Delta_j(\cA)}$, $1\leq j \leq n$. \\
$(ii)$ Let $\cA$ be a strongly hyperbolic monic pencil of degree $n$.
Then $\cA$ is decomposable as 
\begin{equation}
\cA(z)=(z\cI_m-\cY_n)(z\cI_m-\cY_{n-1})\cdots
(z\cI_m-\cY_1),
\end{equation}
with $\spec(\cY_j)\subset \ol {\Delta_j(\cA)}$, $1\leq j \leq n$.
Moreover, each $\cY_j$ is similar to a spectral root $\cZ_j$ of
$\cA$ and hence,
\begin{equation}
\spec(\cY_j)=\spec(\cZ_j)=\spec(\cA)\cap\ol {\Delta_j(\cA)}, 
\quad 1\leq j\leq n. 
\end{equation}
\end{theorem} 

\section{A Class of Matrix-Valued Schr\"odinger Operators \\ with 
Prescribed Finite-Band Spectra} \lb{s4}

This section is devoted to the construction of a class of
matrix-valued Schr\"odinger operators with a prescribed
finite-band spectrum of uniform maximum multiplicity, the principal 
result of this paper.

To begin our analysis we start with a useful result on (scalar) 
Herglotz functions. Even though the result is probably well-known to
experts, we provide an elementary proof for completeness.

Let 
\begin{equation}
\{E_\ell\}_{0\leq \ell\leq 2n}\subseteq \bbR, \; n\in\bbN, \text{ with
$E_\ell<E_{\ell+1}$, $0\leq \ell\leq 2n-1$,} \lb{4.1}
\end{equation}
and introduce the polynomial 
\begin{equation}
R_{2n+1}(z)=\prod_{\ell=0}^{2n} (z-E_\ell), \quad z\in\bbC. \lb{4.2}
\end{equation}
Moreover, we define the square root of $R_{2n+1}$ by 
\begin{equation}
R_{2n+1}(\lambda)^{1/2}=\lim_{\varepsilon\downarrow 0}
R_{2n+1}(\lambda+i\varepsilon)^{1/2}, 
\quad \lambda\in\bbR, \lb{4.3}
\end{equation}
and
\begin{align}
R_{2n+1} (\lambda)^{1/2} &= |R_{2n+1} (\lambda)^{1/2} |\begin{cases}
(-1)^n i& \text{for $\lambda \in (-\infty, E_0)$},\\
(-1)^{n+j} i & \text{for $\lambda \in (E_{2j-1}, E_{2j})$},\; j=1,\dots,
n,\\ (-1)^{n+j}& \text{for $\lambda \in (E_{2j}, E_{2j+1})$}, \;
j=0,\dots, n-1,\\ 1 & \text{for $\lambda \in (E_{2n}, \infty)$},
\end{cases} \no \\
& \hspace*{8.5cm} \lambda\in\bbR \lb{4.4}
\end{align}
and analytically continue $R_{2n+1}^{1/2}$ {}from $\bbR$ to all of
$\bbC\backslash\Sigma$, where $\Sigma$ is defined by
\begin{equation}
\Sigma=\Bigg\{\bigcup_{j=0}^{n-1} [E_{2j},E_{2j+1}]\Bigg\}\cup
[E_{2n},\infty).  \lb{4.4a}
\end{equation}
In this context we also mention the useful formula
\begin{equation}
\ol{R_{2n+1}(\ol z)^{1/2}}=-R_{2n+1}(z)^{1/2}, \quad z\in\bbC_+.
\lb{4.4b} 
\end{equation}

\begin{theorem} \lb{t4.1}
Let $z\in\bbC\backslash\Sigma$ and $n\in\bbN$. Define $R_{2n+1}^{1/2}$ as
in \eqref{4.1}--\eqref{4.4} followed by an analytic continuation to
$\bbC\backslash\Sigma$. Moreover let $F_n$ and $H_{n+1}$ be two monic 
polynomials of degree $n$ and $n+1$, respectively. Then 
\begin{equation}
\f{iF_n(z)}{R_{2n+1}(z)^{1/2}} \lb{4.5}
\end{equation}
is a Herglotz function if and only if all zeros of $F_n$ are real and
there is precisely one zero in each of the intervals $[E_{2j-1},E_{2j}]$,
$1\leq j\leq n$. Moreover, if $iR_{2n+1}^{-1/2}F_n$ is a Herglotz
function, then it can be represented in the form
\begin{equation}
\f{iF_n(z)}{R_{2n+1}(z)^{1/2}}=\f{1}{\pi} \int_{\Sigma} 
\f{F_n(\lambda)d\lambda}{R_{2n+1}(\lambda)^{1/2}}\f{1}{\lambda-z}, 
\quad z\in\bbC\backslash\Sigma. \lb{4.5a}
\end{equation}
Similarly,
\begin{equation}
\f{iH_{n+1}(z)}{R_{2n+1}(z)^{1/2}} \lb{4.6}
\end{equation}
is a Herglotz function if and only if all zeros of $H_{n+1}$ are real and
there is precisely one zero in each of the intervals
$(-\infty,E_0]$ and $[E_{2j-1},E_{2j}]$, $1\leq j\leq n$. Moreover, 
if $iR_{2n+1}^{-1/2}H_{n+1}$ is a Herglotz
function, then it can be represented in the form
\begin{align}
\f{iH_{n+1}(z)}{R_{2n+1}(z)^{1/2}}&=
\Re\bigg(\f{iH_{n+1}(i)}{R_{2n+1}(i)^{1/2}}\bigg) \no \\
& \quad + \f{1}{\pi} \int_{\Sigma} 
\f{H_{n+1}(\lambda)d\lambda}{R_{2n+1}(\lambda)^{1/2}}
\bigg(\f{1}{\lambda-z}-\f{\lambda}{1+\lambda^2}\bigg), 
\quad z\in\bbC\backslash\Sigma. \lb{4.6a}
\end{align}
\end{theorem}
\begin{proof}
We start with the case of $F_n(z)/R_{2n+1}(z)^{1/2}$ in \eqref{4.5}.
Consider a closed counterclockwise oriented contour
$\Gamma_{R,\varepsilon}$ which consists of the semicircle
$C_{\varepsilon}=\{z\in\bbC\,|\, z=E_0+\varepsilon \exp(i\alpha), \,
-\pi/2\leq\alpha\leq \pi/2\}$ centered at $E_0$, the straight line
$L_+=\{z\in\bbC_+\,|\, z=x+i\varepsilon, \, E_0\leq x \leq R\}$, the
following part of the circle of radius $(R^2+\varepsilon^2)^{1/2}$
centered at $E_0$,
$C_R=\{z\in\bbC\,|\, z=E_0+(R^2+\varepsilon^2)^{1/2}\exp(i\beta), \,
\arctan(\varepsilon/R) \leq \beta \leq 2\pi-\arctan(\varepsilon/R)\}$,
and the straight line $L_-=\{z\in\bbC_-\,|\, z=x-i\varepsilon, \, E_0\leq
x\leq R\}$. Then, for $\varepsilon>0$ small enough and $R>0$ sufficiently
large, one infers
\begin{align}
\f{iF_n(z)}{R_{2n+1}(z)^{1/2}}&=\f{1}{2\pi
i}\oint_{\Gamma_{R,\varepsilon}}
\f{1}{\zeta-z}\f{iF_n(\zeta)}{R_{2n+1}(\zeta)^{1/2}}d\zeta \no \\
&\hspace*{-.45cm} \underset{\varepsilon\downarrow 0, R\uparrow\infty}{=}
\f{1}{\pi} \int_{\Sigma} 
\f{1}{\lambda-z}\f{F_n(\lambda)d\lambda}{R_{2n+1}(\lambda)^{1/2}}.
\lb{4.7}
\end{align} 
Here we used \eqref{4.4} to compute the contributions of the contour
integral along $[E_0,R]$ in the limit $\varepsilon\downarrow 0$ and note
that the integral over $C_R$ tends to zero as $R\uparrow 0$ since 
\begin{equation}
\f{F_n(\zeta)}{R_{2n+1}(\zeta)^{1/2}}\underset{\zeta\to\infty}{=}
\Oh\big(|\zeta|^{-1/2}\big). \lb{4.8}
\end{equation}
Next, utilizing the fact that $F_n$ is monic and using \eqref{4.4} again,
one infers that $F_n(\lambda)d\lambda/R_{2n+1}(\lambda)^{1/2}$ represents
a positive measure supported on $\Sigma$ if and only if $F_n$ has
precisely one zero in each of the intervals $[E_{2j-1},E_{2j}]$, $1\leq
j\leq n$. In other words, 
\begin{equation}
\f{F_n(\lambda)}{R_{2n+1}(\lambda)^{1/2}}\geq 0 \text{ on } \Sigma 
\lb{4.9}
\end{equation}
if and only if
$F_n$ has precisely one zero in each of the intervals 
$[E_{2j-1},E_{2j}]$, $1\leq j\leq n$. The Herglotz representation
theorem, Theorem \ref{t2.7}, then finishes the proof of \eqref{4.5a}. The
proof of
\eqref{4.6a} follows along similar  lines taking into account the
additional residues at $\pm i$ inside 
$\Gamma_{R,\varepsilon}$ which are responsible for the real part on
the right-hand side of \eqref{4.6a}.
\end{proof}

Theorem \ref{t4.1} can be improved by invoking ideas developed in
the Appendix of \cite{KN77} (cf.\ also \cite{SY95}). We will pursue this
further in \cite{BGMS02}.

\begin{corollary} \lb{c4.2} 
Let $z\in\bbC\backslash\Sigma$ and $m,n\in\bbN$. Define $R_{2n+1}^{1/2}$
as in
\eqref{4.1}--\eqref{4.4} followed by an analytic continuation to
$\bbC\backslash\Sigma$. Moreover let $\cF_n$ and
$\cH_{n+1}$ be two monic $m\times m$ matrix pencils of degree $n$ and
$n+1$, respectively. Then $(i/2)R_{2n+1}^{-1/2}\cF_n$ is a Herglotz
matrix if and only if the root zones $\Delta_j(\cF_n)$ of $\cF_n$
satisfy
\begin{equation}
\Delta_j(\cF_n)\subseteq [E_{2j-1},E_{2j}], \quad 1\leq j \leq n. 
\lb{4.12a}
\end{equation}
Analogously, $(i/2)R_{2n+1}^{-1/2}\cH_{n+1}$ is a Herglotz matrix if
and only if the root zones $\Delta_j(\cH_{n+1})$ of $\cH_{n+1}$
satisfy
\begin{equation}
\Delta_0(\cH_{n+1})\subset (-\infty,E_{0}], \quad
\Delta_j(\cH_{n+1})\subseteq [E_{2j-1},E_{2j}], \quad 1\leq j \leq n. 
\lb{4.12b}
\end{equation}
If \eqref{4.12a} $($resp., \eqref{4.12b}$)$ holds, then $\cF_n$
$($resp., $\cH_{n+1}$$)$ is a strongly hyperbolic pencil.
\end{corollary}
\begin{proof}
We recall that $(i/2)R_{2n+1}^{-1/2}\cF_n$ is an $m \times m$ Herglotz
matrix if and only if $(f,(i/2)R_{2n+1}^{-1/2}\cF_n f)_{\bbC^m}$ is a
Herglotz function for all $f\in\bbC^m\backslash\{0\}$. Thus, it suffices
to apply Theorem \ref{t4.1}, identifying $F_n$ and $(f,\cF_n 
f)_{\bbC^m}$, to arrive at \eqref{4.12a}. The same argument applied to
$\cH_{n+1}$ yields \eqref{4.12b}. 
\end{proof}

Next, we define the following $2m\times 2m$ matrix $\cM_\Sigma(z,x_0)$
which will turn out to be the underlying Weyl--Titchmarsh matrix
associated with a class of $m \times m$ matrix-valued Schr\"odinger
operators with prescribed finite-band spectra. We introduce, for fixed
$x_0\in\bbR$,
\begin{align}
\cM_\Sigma(z,x_0)&=\big(\cM_{\Sigma,p,q}(z,x_0)
\big)_{1\leq p,q\leq 2} \lb{4.13} \\
&=\f{i}{2R_{2n+1}(z)^{1/2}}\begin{pmatrix} \cH_{n+1,\Sigma}(z,x_0)
& -\cG_{2,n-1,\Sigma}(z,x_0)  \\ 
-\cG_{1,n-1,\Sigma}(z,x_0) & \cF_{n,\Sigma}(z,x_0) \end{pmatrix}, 
\quad z\in\bbC\backslash\Sigma. \no
\end{align}
Here $R_{2n+1}(z)^{1/2}$ is defined as in \eqref{4.1}--\eqref{4.4}
followed by analytic continuation into $\bbC\backslash\Sigma$  and the
polynomial matrix pencils $\cF_{n,\Sigma}$, $\cG_{1,n-1,\Sigma}$,
$\cG_{2,n-1,\Sigma}$, and $\cH_{n+1,\Sigma}$ are introduced as follows: \\

\noindent $(i)$ $\cF_{n,\Sigma}(\cdot,x_0)$ is an $m\times m$ monic
matrix pencil of degree $n$, that is, $\cF_{n,\Sigma}(\cdot,x_0)$
is of the type 
\begin{equation}
\cF_{n,\Sigma}(z,x_0)=\sum_{\ell=0}^n \cF_{n-\ell,\Sigma}(x_0) z^\ell,
\quad \cF_{0,\Sigma}(x_0)=\cI_m, \; z\in\bbC \lb{4.19}
\end{equation}
and 
\begin{equation}
\f{i}{2R_{2n+1}^{1/2}}\cF_{n,\Sigma}(\cdot,x_0) \text{ is assumed to be an
$m\times m$ Herglotz matrix.}  \lb{4.20}  
\end{equation}
Hence $\cF_{n,\Sigma}(\cdot,x_0)$ is a self-adjoint (in fact, strongly
hyperbolic) pencil,
\begin{equation}
\cF_{n,\Sigma}(\ol z,x_0)^*=\cF_{n,\Sigma}(z,x_0), \quad z\in\bbC
\lb{4.19a}
\end{equation}
and $(i/2)R_{2n+1}^{-1/2}\cF_{n,\Sigma}$ and
$2iR_{2n+1}^{1/2}\cF_{n,\Sigma}^{-1}$ admit the Herglotz representations
\begin{align}
&\f{i}{2R_{2n+1}(z)^{1/2}}\cF_{n,\Sigma}(z,x_0)=\f{1}{2\pi} \int_{\Sigma} 
\f{d\lambda}{R_{2n+1}(\lambda)^{1/2}}\cF_{n,\Sigma}(\lambda,x_0) 
\f{1}{\lambda-z}, 
\quad z\in\bbC\backslash\Sigma, \lb{4.21} \\
&iR_{2n+1}(z)^{1/2}\cF_{n,\Sigma}(z,x_0)^{-1} \no \\
&=\f{1}{\pi} \int_\Sigma d\lambda\,
R_{2n+1}(\lambda)^{1/2}\cF_{n,\Sigma}(\lambda,x_0)^{-1} 
\bigg(\f{1}{\lambda -z}-\f{\lambda}{1+\lambda^2}\bigg) \no \\
&\quad +\Gamma_{\Sigma,0}(x_0)-\sum_{k=1}^{N} 
(z-\mu_k (x_0))^{-1} \Gamma_{\Sigma,k}(x_0) , \lb{4.22} \\
& \hspace*{2.2cm} z\in\bbC\backslash
\{\Sigma\cup\{\mu_k(x_0)\}_{1\leq k\leq N}\}, \no
\end{align}
where
\begin{align}
&\Gamma_{\Sigma,0}(x_0)=\Gamma_{\Sigma,0}(x_0)^*\in\bbC^{m\times m}, \quad
0\leq \Gamma_{\Sigma,k}(x_0)\in\bbC^{m\times m}, \; 1\leq k \leq N, \no \\
& \sum_{k=1}^N \rank(\Gamma_{\Sigma,k}(x_0)) \leq mn, \quad 
\mu_k(x_0)\in \bigcup_{j=1}^n [E_{2j-1},E_{2j}], \quad 1\leq k \leq N.
\lb{4.26}
\end{align}
In fact, there are precisely $m$ numbers $\mu_k(x_0)$ in 
$[E_{2j-1},E_{2j}]$ for each $1 \leq j\leq n$, counting multiplicity (they
are the points $z$  where $\cF_{n,\Sigma}(z,x_0)$ is not invertible). \\

\noindent $(ii)$ Given these facts we now define  
\begin{align}
\cG_{1,n-1,\Sigma}(z,x_0)&=\bigg(\sum_{k=1}^N
\f{\varepsilon_k (x_0)}{z-\mu_k(x_0)}
\Gamma_{\Sigma,k}(x_0)\bigg)\cF_{n,\Sigma}(z,x_0),
\lb{4.27} \\
\cG_{2,n-1,\Sigma}(z,x_0)&=\cF_{n,\Sigma}(z,x_0) \bigg(\sum_{k=1}^N
\f{\varepsilon_k (x_0)}{z-\mu_k(x_0)} \Gamma_{\Sigma,k}(x_0)\bigg),
\lb{4.28} \\ 
& \hspace*{-1.95cm} \varepsilon_k (x_0)\in\{1,-1\}, \; 1\leq k \leq N, 
\quad z\in\bbC\backslash\{\mu_k(x_0)\}_{1\leq k\leq N}, \lb{4.29} 
\end{align}
and
\begin{align}
&\cH_{n+1,\Sigma}(z,x_0)=R_{2n+1}(z) \cF_{n,\Sigma}(z,x_0)^{-1}
\lb{4.30}  \\ 
& +\bigg(\sum_{k=1}^N
\f{\varepsilon_k (x_0)}{z-\mu_k(x_0)}
\Gamma_{\Sigma,k}(x_0)\bigg)\cF_{n,\Sigma}(z,x_0) 
\bigg(\sum_{\ell=1}^N
\f{\varepsilon_\ell (x_0)}{z-\mu_\ell(x_0)}
\Gamma_{\Sigma,\ell}(x_0)\bigg),  \no \\
& \hspace*{6.95cm} z\in\bbC\backslash\{\mu_k(x_0)\}_{1\leq k\leq N}. \no
\end{align}

\begin{lemma} \lb{l4.3} Let $z\in\bbC\backslash\{\mu_k(x_0)\}_{1\leq
k\leq N}$. $\cG_{p,n-1,\Sigma}(\cdot,x_0)$, $p=1,2$, are $m\times m$
polynomial matrix pencils of equal degree at most $n-1$ and
$\cH_{n+1,\Sigma}(\cdot,x_0)$ is a self-adjoint $m\times m$ monic
matrix pencil of degree $n+1$. Moreover, the following identities hold.
\begin{align}
&\cG_{2,n-1,\Sigma}(\ol z,x_0)^*=\cG_{1,n-1,\Sigma}(z,x_0), \lb{4.30a}
\\
&\cF_{n,\Sigma}(z,x_0)\cG_{1,n-1,\Sigma}(z,x_0)
=\cG_{2,n-1,\Sigma}(z,x_0)\cF_{n,\Sigma}(z,x_0), \lb{4.31} \\
&\cH_{n+1,\Sigma}(z,x_0)\cG_{2,n-1,\Sigma}(z,x_0)
=\cG_{1,n-1,\Sigma}(z,x_0)\cH_{n+1,\Sigma}(z,x_0), \lb{4.31a} \\
&\cF_{n,\Sigma}(z,x_0)\cH_{n+1,\Sigma}(z,x_0)
-\cG_{2,n-1,\Sigma}(z,x_0)^2
=R_{2n+1}(z)\cI_m, \lb{4.32} \\ 
&\cH_{n+1,\Sigma}(z,x_0)\cF_{n,\Sigma}(z,x_0)
-\cG_{1,n-1,\Sigma}(z,x_0)^2 =R_{2n+1}(z)\cI_m. \lb{4.33} 
\end{align}
\end{lemma}
\begin{proof}
The identities \eqref{4.31}--\eqref{4.33} are obvious {}from
\eqref{4.27}--\eqref{4.30}. Similarly, \eqref{4.30a} is clear {}from
\eqref{4.27}--\eqref{4.29} and \eqref{4.19a}. By \eqref{4.22}, one infers
\begin{equation}
\lim_{z\to \mu_k(x_0)}
\Gamma_{\Sigma,k}(x_0)\cF_{n,\Sigma}(z,x_0)=0=\lim_{z\to
\mu_k(x_0)}\cF_{n,\Sigma}(z,x_0)\Gamma_{\Sigma,k}(x_0) \lb{4.34}
\end{equation}
and hence $\cG_{p,n-1,\Sigma}(\cdot,x_0)$, $p=1,2$, are polynomial
matrix pencils of degree at most $n-1$. By \eqref{4.30a},  
\begin{equation}
\deg(\cG_{1,n-1,\Sigma}(\cdot,x_0))
=\deg(\cG_{2,n-1,\Sigma}(\cdot,x_0))\leq n-1. \lb{4.35}
\end{equation}
Next, using \eqref{4.22} again, one notes that 
\begin{align}
&iR_{2n+1}(\mu_k(x_0))^{1/2}\Gamma_{\Sigma,k}(x_0)=-\Gamma_{\Sigma,k}(x_0)
[(d/dz)\cF_{n,\Sigma}(\mu_k(x_0))]\Gamma_{\Sigma,k}(x_0), \lb{4.38} \\
& \hspace*{9.1cm} 1\leq k\leq N \no
\end{align}
and thus, combining \eqref{4.30} and \eqref{4.38},
\begin{align}
&\res_{z=\mu_k(x_0)} \cH_{n+1,\Sigma}(z,x_0) \lb{4.39} \\
&=iR_{2n+1}(\mu_k(x_0))^{1/2}\Gamma_{\Sigma,k}(x_0)+
\Gamma_{\Sigma,k}(x_0)[(d/dz)\cF_{n,\Sigma}(\mu_k(x_0))]\Gamma_{\Sigma,k}(x_0) 
=0, \no \\ 
& \hspace*{9.9cm} 1\leq k\leq N. \no
\end{align}
Hence, $\cH_{n+1,\Sigma}$ is indeed a polynomial matrix pencil of degree
$n+1$. 
\end{proof}
In fact, $\cH_{n+1,\Sigma}$ is a strongly hyperbolic pencil as shown in
Theorem \ref{t4.9}.

In the following it will be convenient to use the following set of
assumptions.

\begin{hypothesis} \lb{h4.4}
Let $m,n\in\bbN$. Define $R_{2n+1}$ as in \eqref{4.1}, \eqref{4.2} and
$R_{2n+1}^{1/2}$ as in \eqref{4.3}, \eqref{4.4} followed by an analytic
continuation to $\bbC\backslash\Sigma$, with $\Sigma$ introduced in
\eqref{4.4a}. Moreover, let the polynomial $m\times m$ matrix pencils
$\cF_{n,\Sigma}(\cdot,x_0)$, $\cG_{1,n-1,\Sigma}(\cdot,x_0)$,
$\cG_{2,n-1,\Sigma}(\cdot,x_0)$, and $\cH_{n+1,\Sigma}(\cdot,x_0)$ be
defined as in \eqref{4.19}, \eqref{4.20}, \eqref{4.27}--\eqref{4.30}.
\end{hypothesis}

Next, we introduce
\begin{subequations} \lb{4.40}
\begin{align}
&\cM_{\pm,\Sigma}(z,x_0) \no \\
&= \pm
iR_{2n+1}(z)^{1/2}\cF_{n,\Sigma}(z,x_0)^{-1}
-\cG_{1,n-1,\Sigma}(z,x_0)\cF_{n,\Sigma}(z,x_0)^{-1} \lb{4.40a} \\
&=\pm
iR_{2n+1}(z)^{1/2}\cF_{n,\Sigma}(z,x_0)^{-1}
-\cF_{n,\Sigma}(z,x_0)^{-1}\cG_{2,n-1,\Sigma}(z,x_0), 
\lb{4.40b} \\
& \hspace*{5.2cm} 
z\in\bbC\backslash\{\Sigma\cup\{\mu_k(x_0)\}_{1\leq k\leq N}\} \no 
\end{align}
\end{subequations}
and 
\begin{equation}
\cN_{\pm,\Sigma} (z,x_0)=\cM_{-,\Sigma}(z,x_0)\pm \cM_{+,\Sigma}(z,x_0),
\quad z\in\bbC\backslash\{\Sigma\cup\{\mu_k(x_0)\}_{1\leq k\leq N}\}. 
\lb{4.41} 
\end{equation}
We also introduce the open interior $\Sigma^o$ of $\Sigma$
defined by $\Sigma^o=\bigcup_{j=0}^{n-1} (E_{2j},E_{2j+1})\cup
(E_{2n},\infty)$. Then one verifies the following fundamental facts.

\begin{theorem} \lb{t4.5}
Assume Hypothesis \ref{h4.4} and let $z\in\bbC\backslash\{\Sigma\cup 
\{\mu_k(x_0)\}_{1\leq k\leq N}\}$.
Moreover, introduce the $2m\times 2m$ matrix
$\cM_\Sigma(\cdot,x_0)$ as in \eqref{4.13}, \eqref{4.41} and the 
$m\times m$ matrices $\cM_{\pm,\Sigma}(\cdot,x_0)$ as in \eqref{4.40}.
Then, \\
$(i)$ $\pm\cM_{\pm,\Sigma}(\cdot,x_0)$ are $m\times m$ Herglotz
matrices with representations
\begin{align}
\pm \cM_{\pm,\Sigma}(z,x_0)&= \f{1}{\pi} \int_\Sigma d\lambda\,
R_{2n+1}(\lambda)^{1/2}\cF_{n,\Sigma}(\lambda,x_0)^{-1} 
\bigg(\f{1}{\lambda -z}-\f{\lambda}{1+\lambda^2}\bigg) \no \\
& \quad +\Gamma_{\Sigma,0}(x_0) - \sum_{k=1}^N
\f{1 \pm\varepsilon_k (x_0)}{z-\mu_k(x_0)} \Gamma_{\Sigma,k}(x_0),
\lb{4.43} \\ 
& \hspace*{1.75cm} z\in\bbC\backslash\{\Sigma\cup
\{\mu_k(x_0)\}_{1\leq k\leq N}\}. \no 
\end{align}
Moreover, $\cM_{\pm,\Sigma}(\cdot,x_0)$ are the half-line $\cM$-matrices
associated with self-adjoint Schr\"odinger operators
$H^D_{\pm,x_0,\Sigma}$ in $L^2([x_0,\pm\infty))^m$, with a Dirichlet
boundary condition at the point $x_0$ and an $m\times m$ matrix-valued
potential
$\cQ_\Sigma$ satisfying 
\begin{equation}
\cQ_\Sigma=\cQ_\Sigma^*\in L^1_{\loc}(\bbR)^{m\times m}\cap
C^\infty(\bbR\backslash\{x_0\})^{m\times m}, \lb{4.43a}
\end{equation}
given by 
\begin{align}
&H_{\pm,x_0,\Sigma}^D=-\cI_m\f{d^2}{dx^2} + \cQ_\Sigma, \lb{4.44} \\
&\dom(H_{\pm,x_0,\Sigma}^D)=\{g\in L^2((x_0,\pm\infty))^m \,|\,g,g'\in
AC([x_0,c])^m \text{ for all } c\gtrless x_0;  \no \\
& \hspace*{3.45cm} \lim_{\varepsilon\downarrow 0} g(x_0\pm\varepsilon)=0;
\;(-g''+\cQ_\Sigma g)\in L^2((x_0,\pm\infty))^m \}. \no
\end{align}
$(ii)$ The differential expression $\cL_\Sigma=-\cI_m\f{d^2}{dx^2} +
\cQ_\Sigma$ is in the limit point case at $\pm\infty$. \\
$(iii)$ The matrix $\cM_{\Sigma}(\cdot,x_0)$, defined by 
\begin{align}
\cM_\Sigma(z,x_0)&=\big(\cM_{\Sigma,p,q}(z,x_0)\big)_{1\leq p,q\leq 2}
\lb{4.42} 
\\ &=\left(\begin{smallmatrix}
\cM_{\pm,\Sigma}(z,x_0)\cN_{-,\Sigma}(z,x_0)^{-1}\cM_{\mp,\Sigma}(z,x_0)   
\;\; &\cN_{-,\Sigma}(z,x_0)^{-1}\cN_{+,\Sigma}(z,x_0)/2  \\ 
\cN_{+,\Sigma}(z,x_0)\cN_{-,\Sigma}(z,x_0)^{-1}/2  &
\cN_{-,\Sigma}(z,x_0)^{-1} \end{smallmatrix}\right), \no
\end{align} 
is a $2m\times 2m$ Herglotz matrix admitting a representation of the type
\eqref{2.42}, with measure $\Omega_\Sigma(\cdot,x)$ given by 
\begin{equation} 
d\Omega_\Sigma(\lambda,x_0)=\begin{cases} \f{1}{2\pi
R_{2n+1}(\lambda)^{1/2}}\left(\begin{smallmatrix}
\cH_{n+1,\Sigma}(\lambda,x_0) & -\cG_{2,n-1,\Sigma}(\lambda,x_0) \\
-\cG_{1,n-1,\Sigma}(\lambda,x_0) & \cF_{n,\Sigma}(\lambda,x_0) 
\end{smallmatrix}\right)d\lambda, & \lambda\in\Sigma^o, \\
0, & \lambda\in\bbR\backslash\Sigma. \end{cases} \lb{4.42a}
\end{equation}
In addition, $\cM_{\Sigma}(\cdot,x_0)$ is the Weyl--Titchmarsh
$\cM$-matrix associated with the self-adjoint Schr\"odinger operator
$H_\Sigma$ in $L^2(\bbR)^m$ defined by
\begin{align}
&H_\Sigma=-\cI_m \f{d^2}{dx^2}+\cQ_\Sigma, \lb{4.45} \\
&\dom(H_\Sigma)=\{g\in L^2(\bbR)^m \mid g,g^\prime\in
\AC_{\loc}(\bbR)^m;\,
 (-g^{\prime\prime}+\cQ_\Sigma g)\in L^2(\bbR)^m\}. \no
\end{align}
$(iv)$ $H_\Sigma$ has purely absolutely continuous spectrum $\Sigma$,
\begin{equation}
\spec(H_\Sigma)=\spec_{\text{ac}}(H_\Sigma)=\Sigma, \quad 
\spec_{\p}(H_{\Sigma_n})=\spec_{\singc}(H_{\Sigma_n})=\emptyset, \lb{4.46}
\end{equation}
with $\spec(H_{\Sigma_n})$ of uniform spectral multiplicity $2m$.
\end{theorem}
\begin{proof}
The representations \eqref{4.43} for $\cM_{\pm,\Sigma}$ immediately
follow  {}from combining \eqref{4.22}, \eqref{4.27}, and \eqref{4.40a}.
These representations also prove that
$\pm\cM_{\pm,\Sigma}(\cdot,x_0)$ are $m\times m$ Herglotz matrices (cf.\
Theorem \ref{t2.7}). Combining \eqref{4.43} with Theorem \ref{t2.7a}, 
taking into account Lemma~8.3.2 in \cite{Le87}, then yields the
properties stated for $\cQ_\Sigma$. That $\cL_\Sigma=-\cI_m\f{d^2}{dx^2} +
\cQ_\Sigma$ is in the limit point case at $\pm\infty$ can be proved
in analogy to Wienholtz's proof \cite{Wi58} of a result originally due to
Povzner \cite{Po53}, reproduced as Theorem\ 35 in \cite[p.\ 58]{Gl65}. The
corresponding details will be presented in \cite{CG02a}. Equation
\eqref{4.42} follows {}from \eqref{4.13}, \eqref{4.27}--\eqref{4.30}, and
\eqref{4.40}. Relation \eqref{4.46} follows from the explicit formula
\eqref{4.42a} of the spectral measure. In particular, the support
property $\supp(\Omega_\Sigma)=\Sigma$ of the measure $\Omega_\Sigma$ in 
\eqref{4.42a} proves $\spec(H_{\Sigma})=\Sigma$, etc. The uniform maximum
spectral multiplicity $2m$ then follows {}from the fact that
$\rank(d\Omega_\Sigma/d\lambda)=2m$ on the interior $\Sigma^o$ of
$\Sigma$. 
\end{proof}

At this point we cannot yet infer continuity of $\cQ_\Sigma$ at the
boundary point $x_0$. We will subsequently return to this issue in Theorem
\ref{t4.9}. 

In the following we will apply these facts to our concrete class of
matrix-valued Schr\"odinger operators discussed in Theorem \ref{t4.5}. 
In order to find the corresponding Weyl--Titchmarsh matrices
$\cM_{\pm,\Sigma}(z,x)$, we need some preparations. We denote by 
$\psi_{\pm,\Sigma}(z,x,x_0)$ the Weyl solutions \eqref{2.31} associated with
$\cQ_\Sigma$, that is,
\begin{equation}
\psi_{\pm,\Sigma}(z,x,x_0)=\theta_{\Sigma}(z,x,x_0)
+\phi_{\Sigma}(z,x,x_0)
\cM_{\pm,\Sigma}(z,x_0), \quad  z\in\bbC\backslash\Sigma, \lb{4.52a}
\end{equation}
where, in obvious notation, $\theta_{\Sigma}(z,x,x_0)$,
$\phi_{\Sigma}(z,x,x_0)$ denote the fundamental system \eqref{2.16}
corresponding to $\cQ_\Sigma$. Then straightforward computations of
the right-hand sides of \eqref{2.48}--\eqref{2.51} (taking into account
\eqref{2.94}, \eqref{2.95}, \eqref{4.31}, and \eqref{4.33}) yield 
\begin{align}
\cM_\Sigma(z,x)&=\big(\cM_{\Sigma,p,q}(z,x)\big)_{1\leq p,q\leq 2} 
 \lb{4.53} \\
&=\f{i}{2R_{2n+1}(z)^{1/2}}\begin{pmatrix} \cH_{n+1,\Sigma}(z,x)
& -\cG_{2,n-1,\Sigma}(z,x)  \\ 
-\cG_{1,n-1,\Sigma}(z,x) & \cF_{n,\Sigma}(z,x) \end{pmatrix}, 
\quad z\in\bbC\backslash\Sigma, \no
\end{align}
where we abbreviated
\begin{align}
\cF_{n,\Sigma}(z,x)&=\theta_{\Sigma}(z,x,x_0)
\cF_{n,\Sigma}(z,x_0)\theta_{\Sigma}(\ol z,x,x_0)^* \no \\
& \quad + \phi_{\Sigma}(z,x,x_0)\cH_{n+1,\Sigma}(z,x_0)
\phi_{\Sigma}(\ol z,x,x_0)^* \no \\
& \quad - \phi_{\Sigma}(z,x,x_0)\cG_{1,n-1,\Sigma}(z,x_0)
\theta_{\Sigma}(\ol z,x,x_0)^* \no \\
& \quad - \theta_{\Sigma}(z,x,x_0)\cG_{2,n-1,\Sigma}(z,x_0)
\phi_{\Sigma}(\ol z,x,x_0)^*, \lb{4.54} \\
\cG_{1,n-1,\Sigma}(z,x)&=
-\theta_{\Sigma}'(z,x,x_0)\cF_{n,\Sigma}(z,x_0)
\theta_{\Sigma}(\ol z,x,x_0)^*
\no \\ & \quad - \phi_{\Sigma}'(z,x,x_0)\cH_{n+1,\Sigma}(z,x_0)
\phi_{\Sigma}(\ol z,x,x_0)^* \no \\
& \quad + \phi_{\Sigma}'(z,x,x_0)\cG_{1,n-1,\Sigma}(z,x_0)
\theta_{\Sigma}(\ol z,x,x_0)^* \no \\
& \quad + \theta_{\Sigma}'(z,x,x_0)\cG_{2,n-1,\Sigma}(z,x_0)
\phi_{\Sigma}(\ol z,x,x_0)^*, \lb{4.55} \\
\cG_{2,n-1,\Sigma}(z,x)&=-\theta_{\Sigma}(z,x,x_0)
\cF_{n,\Sigma}(z,x_0)\theta_{\Sigma}'(\ol z,x,x_0)^* \no \\
& \quad - \phi_{\Sigma}(z,x,x_0)\cH_{n+1,\Sigma}(z,x_0)
\phi_{\Sigma}'(\ol z,x,x_0)^* \no \\
& \quad + \phi_{\Sigma}(z,x,x_0)\cG_{1,n-1,\Sigma}(z,x_0)
\theta_{\Sigma}'(\ol z,x,x_0)^* \no \\
& \quad + \theta_{\Sigma}(z,x,x_0)\cG_{2,n-1,\Sigma}(z,x_0)
\phi_{\Sigma}'(\ol z,x,x_0)^*, \lb{4.56} \\
\cH_{n+1,\Sigma}(z,x)&=\theta_{\Sigma}'(z,x,x_0)
\cF_{n,\Sigma}(z,x_0)\theta_{\Sigma}'(\ol z,x,x_0)^* \no \\
& \quad + \phi_{\Sigma}'(z,x,x_0)\cH_{n+1,\Sigma}(z,x_0)
\phi_{\Sigma}'(\ol z,x,x_0)^* \no \\
& \quad - \phi_{\Sigma}'(z,x,x_0)\cG_{1,n-1,\Sigma}(z,x_0)
\theta_{\Sigma}'(\ol z,x,x_0)^* \no \\
& \quad - \theta_{\Sigma}'(z,x,x_0)\cG_{2,n-1,\Sigma}(z,x_0)
\phi_{\Sigma}'(\ol z,x,x_0)^*. \lb{4.57}
\end{align}
Considerations of this type can be found in \cite[Sect.\ 8.2]{Le87} 
in the special scalar case $m=1$ and in the matrix context $m\in\bbN$ in 
\cite[Sect.\ 9.4]{Sa99a}.

Differentiating \eqref{4.54}--\eqref{4.57} with respect to $x$ (taking
into account that $\theta''=(\cQ_\Sigma -z\cI_m)\theta$,
$\phi''=(\cQ_\Sigma -z)\phi$) then yields \eqref{4.58}--\eqref{4.61}
below.  Alternatively, these results directly follow from Lemma
\ref{l2.12} identifying $\gg$ and $(i/2)R_{2n+1}^{-1/2}\cF_n$, $\gg_p$ and
$(i/2)R_{2n+1}^{-1/2}\cG_{p,n-1}$, $p=1,2$, and $\gh$ and
$(i/2)R_{2n+1}^{-1/2}\cH_{n+1}$, respectively. 

\begin{lemma} \lb{l4.7}
Assume Hypothesis \ref{h4.4} and let $(z,x)\in\bbC\times\bbR$. Then 
\begin{align}
\cF_{n,\Sigma}'&=-(\cG_{1,n-1,\Sigma}+\cG_{2,n-1,\Sigma}), \lb{4.58} \\
\cG_{1,n-1,\Sigma}'&=-(Q_\Sigma-z\cI_m)\cF_{n,\Sigma}-\cH_{n+1,\Sigma}
&\lb{4.59} \\
&=(-\cF_{n,\Sigma}''+\cF_{n,\Sigma}\cQ_\Sigma-\cQ_\Sigma\cF_{n,\Sigma})/2,
\lb{4.59a} \\
\cG_{1,n-1,\Sigma}''&=-2(\cQ_\Sigma-z\cI_m)\cF_{n,\Sigma}'
-\cQ'_\Sigma\cF_{n,\Sigma}
+\cG_{1,n-1,\Sigma}\cQ_\Sigma-\cQ_\Sigma\cG_{1,n-1,\Sigma},
\lb{4.59b} \\
\cG_{2,n-1,\Sigma}'&=-\cF_{n,\Sigma}(Q_\Sigma-z\cI_m)-\cH_{n+1,\Sigma}
\lb{4.60} \\
&=(-\cF_{n,\Sigma}''+\cQ_\Sigma\cF_{n,\Sigma}-\cF_{n,\Sigma}\cQ_\Sigma)/2,
\lb{4.60a} \\
\cG_{2,n-1,\Sigma}''&=-2\cF_{n,\Sigma}'(\cQ_\Sigma-z\cI_m)
-\cF_{n,\Sigma}\cQ'_\Sigma+\cQ_\Sigma\cG_{2,n-1,\Sigma}
-\cG_{2,n-1,\Sigma}\cQ_\Sigma, \lb{4.60b} \\ 
\cH_{n+1,\Sigma}'&=-\cG_{1,n-1,\Sigma}(Q_\Sigma-z\cI_m)
-(Q_\Sigma-z\cI_m)\cG_{2,n-1,\Sigma}, \lb{4.61} \\
\cH_{n+1,\Sigma}&=[\cF_{n,\Sigma}''-\cF_{n,\Sigma}(\cQ_\Sigma
-z\cI_m)-(\cQ_\Sigma-z\cI_m)\cF_{n,\Sigma}]/2. \lb{4.61a}
\end{align} 
\end{lemma}

In particular, one also verifies the following facts from
\eqref{2.58}--\eqref{2.64}.

\begin{lemma} \lb{l4.8}
Assume Hypothesis \ref{h4.4} and let $(z,x)\in\bbC\times\bbR$. Then
\begin{align}
&\cF_{n,\Sigma}(\ol z,x)^*=\cF_{n,\Sigma}(z,x), \quad  
\cH_{n+1,\Sigma}(\ol z,x)^*=\cH_{n+1,\Sigma}(z,x), \no \\
&\cG_{2,n-1,\Sigma}(\ol z,x)^*=\cG_{1,n-1,\Sigma}(z,x), \lb{4.62} \\
&\cF_{n,\Sigma}(z,x)\cG_{1,n-1,\Sigma}(z,x)
=\cG_{2,n-1,\Sigma}(z,x)\cF_{n,\Sigma}(z,x), \lb{4.63} \\
&\cH_{n+1,\Sigma}(z,x)\cG_{2,n-1,\Sigma}(z,x)
=\cG_{1,n-1,\Sigma}(z,x)\cH_{n+1,\Sigma}(z,x), \lb{4.64} \\
&\cH_{n+1,\Sigma}(z,x)\cF_{n,\Sigma}(z,x)
-\cG_{1,n-1,\Sigma}(z,x)^2=R_{2n+1}(z)\cI_m, \lb{4.65} \\
&\cF_{n,\Sigma}(z,x)\cH_{n+1,\Sigma}(z,x)
-\cG_{2,n-1,\Sigma}(z,x)^2=R_{2n+1}(z)\cI_m. \lb{4.66}
\end{align}
\end{lemma}
\begin{proof}
\eqref{4.62} is clear from \eqref{4.30a}, the fact that
$\cF_{n,\Sigma}(\cdot,x_0)$ and $\cH_{n+1,\Sigma}(\cdot,x_0)$ are
self-adjoint $m\times m$ matrix pencils, \eqref{4.55}, and \eqref{4.56}.
Similarly, \eqref{4.63}--\eqref{4.66} follow from elementary (but
somewhat tedious) calculations directly from \eqref{4.54}--\eqref{4.57},
invoking \eqref{2.72}--\eqref{2.95} and \eqref{4.31}--\eqref{4.33}
repeatedly. 
\end{proof}

Combining \eqref{2.52}--\eqref{2.56} and \eqref{4.53} then yields
\begin{subequations} \lb{4.67}
\begin{align}
&\cM_{\pm,\Sigma}(z,x) \no \\
&= \pm iR_{2n+1}(z)^{1/2}\cF_{n,\Sigma}(z,x)^{-1}
-\cG_{1,n-1,\Sigma}(z,x)\cF_{n,\Sigma}(z,x)^{-1} \lb{4.67a} \\
&=\pm iR_{2n+1}(z)^{1/2}\cF_{n,\Sigma}(z,x)^{-1}
-\cF_{n,\Sigma}(z,x)^{-1}\cG_{2,n-1,\Sigma}(z,x), \lb{4.67b} \\
& \hspace*{7.4cm} \quad z\in\bbC\backslash\bbR. \no 
\end{align}
\end{subequations}
One observes that for each $x\in\bbR$, $\cM_{+,\Sigma}(\cdot,x)$ is the
analytic continuation of $\cM_{-,\Sigma}(\cdot,x)$ through the set
$\Sigma$, and vice versa, 
\begin{align}
&\lim_{\varepsilon\downarrow 0}\cM_{+,\Sigma}(\lambda+i\varepsilon,x)
=\lim_{\varepsilon \downarrow 0}\cM_{-,\Sigma}(\lambda-i\varepsilon,x),
\lb{4.67c} \\
& \hspace*{7mm} \lambda\in \bigcup_{j=0}^{n-1} (E_{2j},E_{2j+1}) \cup
(E_{2n},\infty), \; x\in\bbR. \no
\end{align}
In other words, for each $x\in\bbR$, $\cM_{+,\Sigma}(\cdot,x)$ and 
$\cM_{-,\Sigma}(\cdot,x)$ are the two branches of an analytic
matrix-valued function $\cM_{\Sigma}(\cdot,x)$ on the two-sheeted Riemann
surface of $R_{2n+1}^{1/2}$. Thus, the corresponding potential
$\cQ_\Sigma$ is {\it reflectionless} in the sense discussed in
\cite{CGHL00}, \cite{GT00}, and \cite{KS88}.

Thus, one obtains the following results.

\begin{theorem} \lb{t4.9}
Assume Hypothesis \ref{h4.4} and let $z\in\bbC\backslash\bbR$ and
$x\in\bbR$. Then \\
$(i)$ $\cM_{\pm,\Sigma}(z,\cdot)$ in \eqref{4.67} satisfy the
matrix-valued Riccati-type equation
\begin{equation}
\cM_{\pm,\Sigma}'(z,x)+\cM_{\pm,\Sigma}(z,x)^2=\cQ_\Sigma(x)-z \cI_m,
\quad x\in\bbR,
\;  z\in\bbC\backslash\bbR. \lb{4.68}
\end{equation}
Moreover, $\cM_{\pm,\Sigma}(z,x)$ in \eqref{4.67} are the
$m\times m$ Weyl--Titchmarsh matrices associated with $H^D_{\pm,x,\Sigma}$
on the half-lines $[x,\pm\infty)$ and thus for each
$x\in\bbR$, $\cM_\Sigma(z,x)$ in \eqref{4.53} is a $2m\times 2m$
Weyl--Titchmarsh matrix associated with $H_\Sigma$ on $\bbR$. In
particular, $\cM_{\Sigma}(\cdot,x_0)$ is a $2m\times 2m$
Herglotz matrix of $H_\Sigma$ admitting a representation of the type
\eqref{2.42}, with measure $\Omega_\Sigma(\cdot,x)$ given by 
\begin{equation} 
d\Omega_\Sigma(\lambda,x)=\begin{cases} \f{1}{2\pi
R_{2n+1}(\lambda)^{1/2}}\left(\begin{smallmatrix}
\cH_{n+1,\Sigma}(\lambda,x) & -\cG_{2,n-1,\Sigma}(\lambda,x) \\
-\cG_{1,n-1,\Sigma}(\lambda,x) & \cF_{n,\Sigma}(\lambda,x) 
\end{smallmatrix}\right)d\lambda, & \lambda\in\Sigma^o, \\
0, & \lambda\in\bbR\backslash\Sigma. \end{cases} \lb{4.69}
\end{equation}
$(ii)$ $\cF_{n,\Sigma}(\cdot,x)$ and $\cH_{n+1,\Sigma}(\cdot,x)$ are
strongly hyperbolic $($and hence self-adjoint$)$ $m\times m$ monic matrix
pencils of degree $n$ and $n+1$, respectively, and
$\cG_{p,n-1,\Sigma}(\cdot,x)$, $p=1,2$, are
$m\times m$ matrix pencils of degree $n-1$. \\
$(iii)$ $\cQ_\Sigma\in C^\infty(\bbR)^{m\times m}$.
\end{theorem}
\begin{proof}
\eqref{4.68} is clear from Lemma \ref{l4.7} and \eqref{4.67}. 
Since $\cQ_\Sigma \in C^\infty(\bbR)^{m\times m}$, the initial value
problems 
\begin{subequations} \lb{4.70}
\begin{align}
&\cM_{\pm}'(z,x)+\cM_\pm(z,x)^2=\cQ_\Sigma(x)-z \cI_m,
\quad x\in\bbR,
\;  z\in\bbC\backslash\bbR, \lb{4.70a} \\
&\cM_\pm(z,x_0)=\cM_{\pm,\Sigma}(z,x_0), \lb{4.70b}
\end{align}
\end{subequations}
with $\cM_{\pm,\Sigma}(z,x_0)$ given by \eqref{4.40}, has a unique
solution. Since at $x=x_0$ this solution coincides with the
Weyl--Titchmarsh $M$-matrix $\cM_{\pm,\Sigma}(z,x)|_{x=x_0}$ in
\eqref{4.67} (using the initial condition \eqref{2.15} in
\eqref{4.54}--\eqref{4.57}), $\cM_{\pm,\Sigma}(z,x)$ represents the
Weyl--Titchmarsh matrices  associated with $H^D_{\pm,\Sigma}$ on the
half-lines $[x,\pm\infty)$. This proves part (i). By the known leading
asymptotic behavior \eqref{2.32a} of $\cM_{\pm,\Sigma}(\cdot,x)$ (valid
for each $x_0\in\bbR$, see \cite{CG01}) and that of the diagonal Green's
matrix $\cG_\Sigma(\cdot,x,x)=\cM_{\Sigma,2,2}(\cdot,x)$ of
$H_\Sigma$ as $|z|\to\infty$, $\cF_{n,\Sigma}(\cdot,x)$ and 
$\cH_{n+1,\Sigma}(\cdot,x)$ are monic matrix pencils of degree $n$ and
$n+1$, respectively, and $\cG_{p,n-1,\Sigma}(\cdot,x)$, $p=1,2$, are
$m\times m$ matrix pencils of degree $n-1$. Since the diagonal blocks of
each Herglotz matrix are also Herglotz matrices, one concludes that 
$(i/2)R_{2n+1}^{-1/2}\cF_{n,\Sigma}$ and 
$(i/2)R_{2n+1}^{-1/2}\cH_{n+1,\Sigma}$ are Herglotz matrices. By 
Corollary \ref{c4.2} this then proves that $\cF_{n,\Sigma}$ and
$\cH_{n+1,\Sigma}$ are strongly hyperbolic pencils and hence item (ii)
holds. As in Theorem \ref{t4.5}, $\cQ_\Sigma\in C^\infty((-\infty,x)\cup
(x,\infty))^{m\times m}$. Since $x\in\bbR$ is arbitrary, this proves
(iii). 
\end{proof}

It should be emphasized that the construction of $\cQ_\Sigma$ along the
lines of Section \ref{s4} in the scalar case $m=1$ is due to Levitan
\cite{Le77} (see also \cite{Le77a}, \cite{Le84}, \cite[Ch.\ 8]{Le87},
\cite{LS88}).

\section{Trace Formulas and Connections with the \\Stationary Matrix 
KdV Hierarchy} \lb{s5}

In this section we introduce the stationary matrix Korteweg--de
Vries (KdV) hierarchy and show that the class of finite-band potentials
$\cQ_\Sigma$ constructed in Section \ref{s4} satisfies some (and hence
infinitely many) equations of the stationary KdV equations. We also
introduce trace formulas for KdV invariants.

In order to extend the recursive approach constructing KdV Lax pairs in
the scalar (Abelian) context to the present matrix-valued (non-Abelian)
setting, we focus on an efficient approach introduced by Dubrovin
\cite{Du75} (in the scalar case $m=1$). Recalling
\eqref{2.52}--\eqref{2.64} and Lemma \ref{l2.12}, we state the following
matrix-version of Dubrovin's generating function approach to higher-order
Lax pairs.

We start by introducing the following hypothesis.

\begin{hypothesis} \lb{h5.1}
Fix $m\in\bbN$, suppose $\cQ=\cQ^*\in C^\infty(\bbR)^{m\times
m}$ and introduce the differential expression
\begin{equation}
\cL=-\cI_m\f{d^2}{dx^2}+\cQ, \quad x\in\bbR. \lb{5.1}
\end{equation} 
Suppose $\cL$ is in the limit point case at $\pm\infty$ and introduce the
corresponding self-adjoint operator $H$ in $L^2(\bbR)^{m}$ by
\begin{align}
&H=-\cI_m \f{d^2}{dx^2}+\cQ, \lb{5.2} \\
&\dom(H)=\{g\in L^2(\bbR)^m \mid g,g^\prime\in
\AC_{\loc}(\bbR)^m;\,
 (-g^{\prime\prime}+\cQ g)\in L^2(\bbR)^m\}. \no
\end{align}
\end{hypothesis}

Given Hypothesis \ref{h5.1}, we introduce the generating operator $P_z$ 
by
\begin{align} 
&P_z=\Big(\gg(z,\cdot)\f{d}{dx}+\gg_2(z,\cdot)\Big)(\cL-z\cI_m)^{-1},
\quad z\in\bbC\backslash\bbR, \lb{5.3} \\
&\dom(P_z)=\bigcup_{E\in\bbC\backslash\{\bbR\cup\{z\}\}}\ker(\cL-E\cI_m),
\no 
\end{align}
where $\ker(\cL-E\cI_m)$ denotes the algebraic nullspace of
$\cL-E\cI_m$ (as opposed to the functional analytic
nullspace in $L^2(\bbR)^m$) and $(\cL-z\cI_m)^{-1}$ acts in the obvious
manner by
\begin{equation}
(\cL-z\cI_m)^{-1}\psi=(E-z)^{-1}\psi, \quad
\psi\in\ker(\cL-E\cI_m). \lb{5.4}
\end{equation}
The precise operator theoretic properties of $P_z$ will be irrelevant
in the following.

\begin{lemma} \lb{l5.2}
Assume Hypothesis \ref{h5.1} and let $z\in\bbC\backslash\bbR$,
$x\in\bbR$. Then,
\begin{align}
&\Big[\Big(\gg(z,\cdot)\f{d}{dx}+\gg_2(z,\cdot)\Big)(\cL-z\cI_m)^{-1},
\cL\Big]\bigg|_{\bigcup_{E\in\bbC\backslash\{\bbR\cup\{z\}\}}
\ker(\cL-E\cI_m)} \no \\
&=-2\gg'(z,\cdot)\bigg|_{\bigcup_{E\in\bbC\backslash\{\bbR\cup\{z\}\}}
\ker(\cL-E\cI_m)}. \lb{5.5}
\end{align}
\end{lemma}
\begin{proof}
Fix $\psi\in\ker(\cL-E\cI_m)$ for some
$E\in\bbC\backslash\{\bbR\cup\{z\}\}$. Then one computes
\begin{align}
&\Big[\Big(\gg(z,\cdot)\f{d}{dx}+\gg_2(z,\cdot)\Big)(\cL-z\cI_m)^{-1},
\cL\Big]\psi \no \\
&=(E-z)^{-1}\{[2\gg_2'+\gg''+\gg\cQ-\cQ\gg]\psi' \no \\
& \quad\;
[\gg_2''+2\gg'(\cQ-z\cI_m)+\gg\cQ'+\gg_2\cQ-\cQ\gg_2+2(z-E)\gg']\psi\} 
\no \\
& =-2\gg'\psi, \lb{5.6}
\end{align}
using \eqref{2.67a} and \eqref{2.67b}. 
\end{proof}

A variant of Dubrovin's idea of a generating operator for KdV Lax pairs
was also used by Olmedilla, Mart\'{\i}nez Alonso, and Guil \cite{OMG81}.
Their approach, however, focuses on formal operator expansions and formal
pseudo-differential operators. For a different approach we refer to
\cite[Ch.~15]{Di91}, \cite{GD77}.

Next, we recall that 
\begin{align}
\gg(z,x)&=[\cM_-(z,x)-\cM_+(z,x)]^{-1}, \lb{5.7} \\
\gg_1(z,x)&=\f12
[\cM_-(z,x)+\cM_+(z,x)][\cM_-(z,x)-\cM_+(z,x)]^{-1}, \lb{5.8b} \\
\gg_2(z,x)&=\f12
[\cM_-(z,x)-\cM_+(z,x)]^{-1}[\cM_-(z,x)+\cM_+(z,x)], \lb{5.8} \\
\gh(z,x)&=\cM_\pm(z,x)[\cM_-(z,x)-\cM_+(z,x)]^{-1}\cM_\mp(z,x), \lb{5.8a} 
\end{align}
and note that by Theorem \ref{t2.10} the right-hand sides of \eqref{5.7}
and \eqref{5.8} admit asymptotic expansions in cones avoiding the
spectrum of $H$. In particular, one thus obtains the asymptotic expansions
\begin{align}
\gg(z,x)&\underset{\substack{\abs{z}\to\infty\\ z\in
C_\varepsilon}}{=}\f{i}{2z^{1/2}}\sum_{k=0}^\infty
\hatt\gR_{k}(x)z^{-k}, \quad \hatt\gR_0(x)=\cI_m, \lb{5.9} \\
\gg_p(z,x)&\underset
{\substack{\abs{z}\to\infty\\ z\in
C_\varepsilon}}{=}\f{i}{2z^{1/2}}\sum_{k=0}^\infty 
\hatt\gG_{p,k}(x)z^{-k}, \quad \hatt\gG_{p,0}(x)=0, \quad p=1,2, \lb{5.10}
\\
\gh(z,x)&\underset{\substack{\abs{z}\to\infty\\ z\in
C_\varepsilon}}{=}\f{iz^{1/2}}{2}\sum_{k=0}^\infty
\hatt\gH_{k}(x)z^{-k}, \quad \hatt\gH_0(x)=\cI_m, \lb{5.10a}
\end{align}
for some coeficients $\hatt\gR_k$, $\hatt\gG_{p,k}$, $p=1,2$, and
$\hatt\gH_k$, which are universal differential polynomials in $\cQ$ by
Remark
\ref{r2.11}\,(i). Explicitly, one obtains
\begin{align}
\hatt\gR_0&=\cI_m, \quad \hatt\gR_1=\tfrac{1}{2}\cQ, \quad 
\hatt\gR_2=-\tfrac{1}{8}\cQ''+\tfrac{3}{8}\cQ^2, \lb{5.10b} \\
\hatt\gG_{1,0}&=-\tfrac{1}{4}\cQ', \quad
\hatt\gG_{1,1}=\tfrac{1}{16}\cQ'''-\tfrac{1}{8}(\cQ^2)'
-\tfrac{1}{8}\cQ'\cQ, \lb{5.10c} \\
\hatt\gG_{2,0}&=-\tfrac{1}{4}\cQ', \quad
\hatt\gG_{2,1}=\tfrac{1}{16}\cQ'''-\tfrac{1}{8}(\cQ^2)'
-\tfrac{1}{8}\cQ\cQ', \lb{5.10d} \\
\hatt\gH_0&=\cI_m, \quad \hatt\gH_1=-\tfrac{1}{2}\cQ, \quad 
\hatt\gH_2=\tfrac{1}{8}\cQ''-\tfrac{1}{8}\cQ^2, \lb{5.10e} \\
& \text{ etc.} \no
\end{align}

Motivated by Lemma \ref{l5.2}, we now introduce the
$m\times m$ matrix-valued differential expressions $\hatt\cP_{2k+1}$ by
\begin{equation}
\hatt\cP_{2k+1}=\sum_{\ell=0}^k\Big(\hatt\gR_\ell\f{d}{dx}
+\hatt\gG_{2,\ell}\Big)\cL^{k-\ell}, \quad k\in\bbN_0. \lb{5.11}
\end{equation}
In analogy to Lemma \ref{l5.2} one then obtains the following result.

\begin{lemma} \lb{l5.2a}
Assume Hypothesis \ref{h5.1} and let $z\in\bbC\backslash\bbR$,
$x\in\bbR$. Then,
\begin{equation}
\big[\hatt\cP_{2k+1},\cL\big]=2\hatt\gR_{k+1}', \quad k\in\bbN_0.
\lb{5.12}
\end{equation}
\end{lemma}
\begin{proof}
Assuming $\psi\in\ker(\cL-z\cI_m)$ one computes
\begin{align}
&\big[\hatt\cP_{2k+1},\cL\big]\psi =\sum_{\ell=0}^k
z^{k-\ell}\big(2\hatt\gG_{2,\ell}'+\hatt\gR_\ell''+\hatt\gR_\ell\cQ
-\cQ\hatt\gR_\ell\big)\psi' \no \\
&+\sum_{\ell=0}^k
z^{k-\ell}\big(\hatt\gG_{2,\ell}''+2\hatt\gR_\ell'(\cQ-z\cI_m)+ 
\hatt\gR_\ell\cQ'+\hatt\gG_{2,\ell}\cQ-\cQ\hatt\gG_{2,\ell}\big)\psi \no
\\ 
&=\sum_{\ell=0}^k
z^{k-\ell}\big(2\hatt\gG_{2,\ell}'+\hatt\gR_\ell''+\hatt\gR_\ell\cQ
-\cQ\hatt\gR_\ell\big)\psi' -2\hatt\gR_0'\psi \no \\
&+\sum_{\ell=0}^k
z^{k-\ell}\big(\hatt\gG_{2,\ell}''-2\hatt\gR_{\ell+1}'+2\hatt\gR_\ell'\cQ
+ \hatt\gR_\ell\cQ'+\hatt\gG_{2,\ell}\cQ-\cQ\hatt\gG_{2,\ell}\big)\psi + 
2\hatt\gR_{k+1}'\psi\no \\ 
&=2\hatt\gR_{k+1}'\psi. \lb{5.13}
\end{align}
Here we used $\hatt\gR_0'=0$ and 
\begin{align}
&2\hatt\gG_{2,\ell}'+\hatt\gR_\ell''+\hatt\gR_\ell\cQ
-\cQ\hatt\gR_\ell=0, \quad \ell\in\bbN_0, \lb{5.13a} \\
&\hatt\gG_{2,\ell}''-2\hatt\gR_{\ell+1}'+2\hatt\gR_\ell'\cQ
+ \hatt\gR_\ell\cQ'+\hatt\gG_{2,\ell}\cQ-\cQ\hatt\gG_{2,\ell}=0, \quad
\ell\in\bbN_0, \lb{5.13b} 
\end{align}
which follow from inserting the asymptotic expansion \eqref{5.9} and
\eqref{5.10} into \eqref{2.67a} and \eqref{2.67b} (which is permitted by
Theorem \ref{t2.10}). Relation \eqref{5.13} implies \eqref{5.12} since
$\hatt\cP_{2k+1}$ and $\cL$ are $m\times m$ matrix-valued differential
expressions of finite-order while
$\bigcup_{z\in\bbC}\ker(\cL-z\cI_m)$ is an infinite-dimensional space of
$C^\infty(\bbR)^m$-functions. 
\end{proof}

Introducing 
\begin{equation}
\cP_{2k+1}=\sum_{\ell=0}^k c_{k-\ell}\hatt\cP_{2\ell+1}, \quad
k\in\bbN_0, \lb{5.14}
\end{equation}
where 
\begin{equation}
\{c_\ell\}_{\ell=1,\dots,k}\subset\bbC, \quad c_0=1 \lb{5.15}
\end{equation}
denotes a set of constants, the pairs $(\cP_{2k+1},\cL)$, $k\in\bbN_0$,
by definition, represent the Lax pairs of the (matrix-valued) KdV
hierarchy. More precisely, varying $\{c_\ell\}_{\ell\in\bbN}\subset\bbC$,
the set of evolution equations,
\begin{equation}
\f{d}{dt}\cL-[\cP_{2k+1},\cL]=0, \quad k\in\bbN_0, \lb{5.16}
\end{equation}
or equivalently, the set of equations,
\begin{equation}
\kdv_{k}(\cQ)=\cQ_t-2\sum_{\ell=0}^k c_{k-\ell}
\hatt\gR_{\ell+1}'(\cQ,\dots)=0, \quad k\in\bbN_0, \lb{5.17}
\end{equation}
represents the (matrix-valued) KdV hierarchy of evolution equations. The
corresponding stationary KdV hierarchy, characterized by 
\begin{equation}
\cQ_t=0, \text{ or equivalently, by } [\cP_{2k+1},\cL]=0, \quad
k\in\bbN_0, \lb{5.17a}
\end{equation}
is then given by
\begin{equation}
\skdv_{k}(\cQ)=-2\sum_{\ell=0}^k c_{k-\ell}
\hatt\gR_{\ell+1}'(\cQ,\dots)=0, \quad k\in\bbN_0. \lb{5.17b}
\end{equation}

\begin{remark} \lb{r5.3}
By Remark \ref{r2.11}, each $\hatt\gR_\ell$ is a universal polynomial in
$\cQ$ and its $x$-derivatives and occasionally we slightly abuse notation
and indicate this by writing $\hatt\gR_\ell(\cQ,\dots)$ for
$\hatt\gR_\ell(x)$, $\hatt\gR_{\ell+1}'(\cQ,\dots)$ for
$\hatt\gR_{\ell+1}'(x)$, etc. Explicit formulas for the differential
polynomials $\hatt\gR_\ell(\cQ,\dots)$ were derived for instance, in
\cite{AS00} and \cite{Po01}.
\end{remark}

In order to make the connection with the finite-band formalism of Section
\ref{s4} we now recall 
\begin{align}
&\f{i}{2R_{2n+1}(z)^{1/2}}\cF_{n,\Sigma}(z,x)=[\cM_{-,\Sigma}(z,x)
-\cM_{+,\Sigma}(z,x)]^{-1}, \lb{5.18} \\
&\f{i}{2R_{2n+1}(z)^{1/2}}\cG_{1,n-1,\Sigma}(z,x) \no \\
&\quad =\f12 [\cM_{-,\Sigma}(z,x)+\cM_{+,\Sigma}(z,x)]
[\cM_{-,\Sigma}(z,x)-\cM_{+,\Sigma}(z,x)]^{-1}, \lb{5.18a} \\
&\f{i}{2R_{2n+1}(z)^{1/2}}\cG_{2,n-1,\Sigma}(z,x) \no \\
&\quad =\f12 [\cM_{-,\Sigma}(z,x)-\cM_{+,\Sigma}(z,x)]^{-1} 
[\cM_{-,\Sigma}(z,x)+\cM_{+,\Sigma}(z,x)], \lb{5.19} \\
&\f{i}{2R_{2n+1}(z)^{1/2}}\cH_{n+1,\Sigma}(z,x) \no \\
&\quad = \cM_{\pm,\Sigma}(z,x)
[\cM_{-,\Sigma}(z,x)-\cM_{+,\Sigma}(z,x)]^{-1}\cM_{\mp,\Sigma}(z,x),
\lb{5.19a} 
\end{align}
and note that \eqref{5.18}--\eqref{5.19a} admit expansions convergent 
in a neighborhood of infinity. In particular, 
\begin{align}
\f{1}{R_{2n+1}(z)^{1/2}}\cF_{n,\Sigma}(z,x)&\underset
{\substack{\abs{z}\to\infty\\ z\in
C_\varepsilon}}{=}\f{1}{z^{1/2}}\sum_{k=0}^\infty
\hatt\gR_{k,\Sigma}(x)z^{-k}, \quad \hatt\gR_{0,\Sigma}(x)=\cI_m,
\lb{5.20} \\
\f{1}{R_{2n+1}(z)^{1/2}}\cG_{p,n-1,\Sigma}(z,x)&\underset
{\substack{\abs{z}\to\infty\\ z\in
C_\varepsilon}}{=}\f{1}{z^{1/2}}\sum_{k=0}^\infty 
\hatt\gG_{p,k,\Sigma}(x)z^{-k}, \quad \hatt\gG_{p,0,\Sigma}(x)=0, 
\; p=1,2, \lb{5.21} \\
\f{1}{R_{2n+1}(z)^{1/2}}\cH_{n+1,\Sigma}(z,x)&\underset
{\substack{\abs{z}\to\infty\\ z\in
C_\varepsilon}}{=}\f{1}{z^{1/2}}\sum_{k=0}^\infty
\hatt\gH_{k,\Sigma}(x)z^{-k}, \quad \hatt\gH_{0,\Sigma}(x)=\cI_m,
\lb{5.21a} 
\end{align}
for $|z|$ sufficiently large. Here the coefficients $\hatt\gR_{k,\Sigma}$
and $\hatt\gG_{2,k,\Sigma}$ are the universal differential polynomials
$\hatt\gR_k=\hatt\gR_k(\cQ_\Sigma,\dots)$ and
$\hatt\gG_{2,k}=\hatt\gG_{2,k}(\cQ_\Sigma,\dots)$ in \eqref{5.9} and
\eqref{5.10} (with $\cQ$ replaced by $\cQ_\Sigma$). We also recall 
\begin{align}
\cF_{n,\Sigma}(z,x)&=\sum_{\ell=0}^n \cF_{n-\ell,\Sigma}(x) z^\ell,
\quad \cF_{0,\Sigma}(x)=\cI_m,  \lb{5.22} \\
\cG_{p,n-1,\Sigma}(z,x)&=\sum_{\ell=0}^{n-1} \cG_{p,n-1-\ell,\Sigma}(x)
z^\ell, \quad p=1,2, \lb{5.23} \\
\cH_{n+1,\Sigma}(z,x)&=\sum_{\ell=0}^n \cH_{n+1-\ell,\Sigma}(x) z^\ell,
\quad \cH_{0,\Sigma}(x)=\cI_m.  \lb{5.23a}
\end{align}
Since we seek the connection between the set of coefficients
$\hatt\gR_{k,\Sigma}$, $\hatt\gG_{p,k,\Sigma}$, $\hatt\gH_{k,\Sigma}$ and
$\cF_{k,\Sigma}$, $\cG_{p,k,\Sigma}$, $\cH_{k,\Sigma}$, we next consider
the following elementary expansions. Let
\begin{equation}
\eta\in\bbC \; \text{ such that } \; 
|\eta|<\min\{|E_0|^{-1},\dots, |E_{2n}|^{-1}\}. \lb{5.24}
\end{equation}
Then 
\begin{equation}
\bigg(\prod_{\ell=0}^{2n} \big(1-E_\ell\eta \big)
\bigg)^{-1/2}=\sum_{k=0}^{\infty}\hat c_k(\ul E)\eta^{k}, \lb{5.25}
\end{equation}
where
\begin{align}
\hat c_0(\ul E)&=1, \no \\
\hat c_k(\ul E)&=\sum_{\substack{j_0,\dots,j_{2n}=0\\
 j_0+\cdots+j_{2n}=k}}^{k}
\f{(2j_0)!\cdots(2j_{2n})!}
{2^{2k} (j_0!)^2\cdots (j_{2n}!)^2}E_0^{j_0}\cdots E_{2n}^{j_{2n}},
\quad k\in\bbN. \lb{5.26}
\end{align} 
The first few coefficients explicitly read
\begin{align}
\hat c_0(\ul E)&=1, \; 
\hat c_1(\ul E)=\f12\sum_{\ell=0}^{2n} E_\ell, \no \\
\hat c_2(\ul E)&=\f14\sum_{\substack{\ell_1,\ell_2=0\\ \ell_1<
\ell_2}}^{2n} E_{\ell_1} E_{\ell_2}+\f38 \sum_{\ell=0}^{2n} E_\ell^2, 
\quad \text{etc.} \lb{5.28}
\end{align}
Similarly, one has 
\begin{equation}
\bigg(\prod_{\ell=0}^{2n} \big(1-E_\ell \eta \big)
\bigg)^{1/2}=\sum_{k=0}^{\infty}c_k(\ul E)\eta^{k}, \lb{5.29}
\end{equation}
where
\begin{align}
c_0(\ul E)&=1,\no \\
c_k(\ul E)&=-\!\!\!\!\!\sum_{\substack{j_0,\dots,j_{2n}=0\\
 j_0+\cdots+j_{2n}=k}}^{k}\!\!\!\!\!
\f{(2j_0)!\cdots(2j_{2n})!}
{2^{2k} (j_0!)^2\cdots (j_{2n}!)^2 (2j_0-1)\cdots(2j_{2n}-1)}
E_0^{j_0}\cdots E_{2n}^{j_{2n}}, \no \\
& \hspace*{9cm} k\in\bbN. \label{5.30} 
\end{align}
The first few coefficients explicitly are given by
\begin{align}
 c_0(\ul E)&=1, \quad 
c_1(\ul E)=-\f12\sum_{\ell=0}^{2n} E_\ell, \no \\ 
c_2(\ul E)&=\f14\!\sum_{\substack{\ell_1,\ell_2=0\\ \ell_1<
\ell_2}}^{2n}\!\!\!\! E_{\ell_1} E_{\ell_2}-\f18 \sum_{\ell=0}^{2n}
E_\ell^2, \quad \text{etc.} \lb{5.31}
\end{align}
\begin{lemma} \lb{l5.4}
Assume Hypothesis \ref{h5.1} and let $x\in\bbR$. Then,
\begin{align}
\cF_{\ell,\Sigma}(x)&=\sum_{k=0}^\ell c_{\ell-k}(\ul E)\hatt
\gR_{k,\Sigma}(x), \quad \ell=0,\dots,n, \lb{5.32} \\
\hatt\gR_{\ell,\Sigma}(x)&=\sum_{k=0}^\ell \hat c_{\ell-k}(\ul
E)\cF_{k,\Sigma}(x), \quad \ell=0,\dots,n, \lb{5.33} \\
\cG_{p,\ell,\Sigma}(x)&=\sum_{k=0}^{\ell} c_{\ell-k}(\ul E)
\hatt \gG_{p,k+1,\Sigma}(x), \quad \ell=0,\dots,n-1, \;\, p=1,2, \lb{5.34}
\\
\hatt\gG_{p,\ell+1,\Sigma}(x)&=\sum_{k=0}^{\ell} \hat c_{\ell-k}(\ul E)
\cG_{p,k,\Sigma}(x), \quad \ell=0,\dots,n-1, \;\, p=1,2, \lb{5.35} \\
\cH_{\ell,\Sigma}(x)&=\sum_{k=0}^\ell c_{\ell-k}(\ul E)\hatt
\gH_{k,\Sigma}(x), \quad \ell=0,\dots,n+1, \lb{5.35a} \\
\hatt\gH_{\ell,\Sigma}(x)&=\sum_{k=0}^\ell \hat c_{\ell-k}(\ul
E)\cH_{k,\Sigma}(x), \quad \ell=0,\dots,n+1. \lb{5.35b} 
\end{align}
\end{lemma}
\begin{proof}
It suffices to prove \eqref{5.32} and \eqref{5.33} and so we omit the
analogous proofs of \eqref{5.34}--\eqref{5.35b}. Since for $|z|$
sufficiently large,
\begin{align}
z^{-n}\cF_{n,\Sigma}(z,x)&=\sum_{\ell=0}^n
\cF_{\ell,\Sigma}(x)z^{-\ell} \no \\
&=z^{-n-(1/2)}R_{2n+1}(z)^{1/2} \sum_{\ell=0}^\infty
\hatt\gR_{\ell,\Sigma}(x)z^{-\ell} \no \\ 
&=\sum_{k=0}^\infty c_k(\ul E)z^{-k}\sum_{\ell=0}^\infty
\hatt\gR_{\ell,\Sigma}(x)z^{-\ell} \no \\
&=\sum_{\ell=0}^\infty \bigg(\sum_{k=0}^\ell c_{\ell-k}(\ul
E) \hatt\gR_{k,\Sigma}(x)\bigg)z^{-\ell} \lb{5.36}
\end{align}
and hence \eqref{5.32}. Equation \eqref{5.33} is then clear from
\eqref{5.32} and 
\begin{equation}
\sum_{\ell=0}^k \hat c_{k-\ell}(\ul E)c_\ell(\ul E) =\delta_{k,0}, 
\quad k\in\bbN_0. \lb{5.37}
\end{equation}
The latter follows from multiplying \eqref{5.25} and \eqref{5.29},
comparing coefficients of $\eta^{k}$.
\end{proof}

Given these preliminaries we can now state the principal result of this
section.

\begin{theorem} \lb{t5.5}
The self-adjoint finite-band potential $\cQ_\Sigma\in
C^\infty(\bbR)^{m\times m}$, discussed in Theorems \ref{t4.5} and
\ref{t4.9}, is a stationary KdV solution satisfying
\begin{equation}
\skdv_{n}(\cQ_\Sigma)=-2\sum_{\ell=0}^n c_{n-\ell}(\ul E)
\hatt\gR_{\ell+1}'(\cQ_\Sigma,\dots)=0, \lb{5.41}
\end{equation}
with $c_\ell(\ul E)$ defined in \eqref{5.30} and $\hatt\gR_{\ell+1}$
the universal differential polynomials $($with respect to $\cQ$$)$ in
\eqref{5.9}.
\end{theorem}  
\begin{proof}
Introducing the $m\times m$ matrix-valued differential expression
$\cP_{2k+1,\Sigma}$ by
\begin{equation}
\cP_{2n+1,\Sigma}=\sum_{\ell=0}^n c_{n-\ell}(\ul E)\hatt\cP_{2\ell+1}
=\sum_{\ell=0}^n \Big(\cF_{n-\ell,\Sigma}(\cdot)\f{d}{dx}+
\cG_{2,n-1-\ell,\Sigma}(\cdot)\Big)\cL^\ell \lb{5.42}
\end{equation}
(cf.\ \eqref{5.32}, \eqref{5.34}, and \eqref{5.14}), one computes for
$\psi\in\ker(\cL-z\cI_m)$, $z\in\bbC$, 
\begin{equation}
\cP_{2n+1,\Sigma}\psi=\cF_{n,\Sigma}(z,x)\psi'+ 
\cG_{2,n-1,\Sigma}(z,x)\psi, \lb{5.43} 
\end{equation}
and hence,
\begin{align}
&[\cP_{2n+1,\Sigma},\cL]\psi=[2\cG_{2,n-1,\Sigma}+\cF_{n,\Sigma}''
\cF_{n,\Sigma}\cQ_\Sigma-\cQ_\Sigma\cF_{n,\Sigma}]\psi' \lb{5.44} \\
&+[\cG_{2,n-1,\Sigma}''+2\cF_{n,\Sigma}'(\cQ_\Sigma-z\cI_m)+
\cF_{n,\Sigma}\cQ_\Sigma'+\cG_{2,n-1,\Sigma}\cQ_\Sigma
-\cQ_\Sigma\cG_{2,n-1,\Sigma}]\psi=0 \no
\end{align}
by \eqref{4.60a} and \eqref{4.60b}. Since $z\in\bbC$ is arbitrary, this
implies
\begin{equation}
[\cP_{2n+1,\Sigma},\cL]=0, \lb{5.45}
\end{equation}
completing the proof by \eqref{5.17a}, \eqref{5.17b}. 
\end{proof}

Next, we turn to a discussion of trace formulas for the finite-band
potential $\cQ_\Sigma$ in terms of (matrix) roots of $\cF_{n,\Sigma}$ 
and $\cH_{n+1,\Sigma}$. 

\begin{theorem} \lb{t5.7}
Let $(z,x)\in\bbC\times\bbR$ and assume $\cQ_\Sigma$ to be the
self-adjoint finite-band potential discussed in Theorems \ref{t4.5} and
\ref{t4.9}. In addition, let the monic self-adjoint matrix pencils
$\cF_{n,\Sigma}$ and
$\cH_{n+1,\Sigma}$ be given by \eqref{4.54} and \eqref{4.57}. Then
$\cF_{n,\Sigma}(\cdot,x)$ and
$\cH_{n+1,\Sigma}(\cdot,x)$ are strongly hyperbolic pencils and hence
admit the factorizations
\begin{align}
\cF_{n,\Sigma}(z,x)&=(z\cI_m-\cU_n(x))(z\cI_m-\cU_{n-1}(x))\cdots
(z\cI_m-\cU_1(x)), \lb{5.58} \\
\cH_{n+1,\Sigma}(z,x)&=(z\cI_m-\cV_n(x))(z\cI_m-\cV_{n-1}(x))\cdots
(z\cI_m-\cV_0(x)), \lb{5.59}
\end{align}
with 
\begin{align}
\spec(\cV_0(x))&=\spec(\cH_{n+1,\Sigma})\cap\ol{\Delta_0(\cH_{n+1,\Sigma})}
\subset (-\infty,E_0], \lb{5.60} \\
\spec(\cU_j(x))&=\spec(\cF_{n,\Sigma})\cap\ol{\Delta_j(\cF_{n,\Sigma})}
\subseteq [E_{2j-1},E_{2j}], \quad 1\leq j\leq n, \lb{5.61} \\
\spec(\cV_j(x))&=\spec(\cH_{n+1,\Sigma})\cap\ol{\Delta_0(\cH_{n+1,\Sigma})}
\subseteq [E_{2j-1},E_{2j}], \quad 1\leq j\leq n. \lb{5.62}
\end{align}
Moreover, one obtains the sequence of trace formulas
\begin{align}
\sum_{\substack{j_1,j_2,\dots,j_k=0\\j_1<j_2<\dots<j_k}}^n
\cU_{j_k}(x)\cdots \cU_{j_2}(x)\cU_{j_1}(x)&=\sum_{\ell=0}^k
c_{k-\ell}(\ul E) \hatt \gR_{\ell,\Sigma}(x), \quad 1\leq k\leq n,
\lb{5.62a} \\ 
0&=\sum_{\ell=0}^k c_{k-\ell}(\ul E) \hatt
\gR_{\ell,\Sigma}(x), \quad k\geq n+1, \lb{5.62b} \\
\sum_{\substack{j_0,j_1,\dots,j_k=0\\j_0<j_1<\dots<j_k}}^n
\cV_{j_k}(x)\cdots \cV_{j_1}(x)\cV_{j_0}(x)&=\sum_{\ell=0}^k
c_{k-\ell}(\ul E)
\hatt \gH_{\ell,\Sigma}(x), \quad 1\leq k\leq n, \lb{5.62c} \\
0&=\sum_{\ell=0}^k c_{k-\ell}(\ul E) \hatt
\gH_{\ell,\Sigma}(x), \quad k\geq n+1. \lb{5.62d}
\end{align}
In particular, in the special case $k=1$, $\cQ_\Sigma$ satisfies the trace
formulas 
\begin{align}
\cQ_\Sigma(x)&=\bigg(\sum_{\ell=0}^{2n}
E_\ell\bigg)\cI_m-2\sum_{j=1}^n\cU_j(x), \lb{5.63} \\
&=-\bigg(\sum_{\ell=0}^{2n}
E_\ell\bigg)\cI_m+2\sum_{k=0}^n\cV_k(x). \lb{5.64} 
\end{align}
In addition, one obtains
\begin{equation}
\cQ_\Sigma^{(r)}\in C^\infty(\bbR)^{m\times m} \cap
L^\infty(\bbR)^{m\times m} \text{ for all $r\in\bbN_0$}. \lb{5.65}
\end{equation}
\end{theorem}
\begin{proof}
By \eqref{5.18} and \eqref{5.19a},
$i/R_{2n+1}^{1/2}\cF_{n,\Sigma}(\cdot,x)$ and
$i/R_{2n+1}^{-1/2}\cH_{n+1,\Sigma}(\cdot,x)$ are Herglotz matrices and
hence \eqref{4.12a} and \eqref{4.12b} apply. In particular,
$\cF_{n,\Sigma}(\cdot,x)$ and $\cH_{n+1,\Sigma}(\cdot,x)$ are strongly
hyperbolic pencils. Since both are monic, the factorizations \eqref{5.58}
and \eqref{5.59}, as well as \eqref{5.60}--\eqref{5.62}, hold by Theorem
\ref{t3.5}\,(ii). By \eqref{5.20} one infers
\begin{align}
z^{-n}\cF_{n,\Sigma}(z,x)=&\sum_{k=0}^n 
\bigg(\sum_{\substack{j_1,j_2,\dots,j_k=0\\
j_1<j_2<\dots<j_k}}^n
\cU_{j_k}(x)\cdots \cU_{j_2}(x)\cU_{j_1}(x)\bigg)z^{-k} \no \\
\underset{\substack{\abs{z}\to\infty\\ z\in
C_\varepsilon}}{=}&z^{-n-(1/2)}R_{2n+1}(z)^{1/2}\sum_{k=0}^\infty
\hatt\gR_{k,\Sigma}(x)z^{-k} \no \\
\underset{\substack{\abs{z}\to\infty\\ z\in
C_\varepsilon}}{=}&\bigg(\prod_{\ell=0}^{2n}
\big(1-(E_\ell/z)\big)\bigg)^{1/2}\sum_{k=0}^\infty \hatt
\gR_{k,\Sigma}(x)z^{-k} \no \\
\underset{\substack{\abs{z}\to\infty\\ z\in
C_\varepsilon}}{=}&\sum_{k=0}^\infty\bigg(\sum_{\ell=0}^k
c_\ell(\ul E)\hatt\gR_{k-\ell,\Sigma}(x)\bigg)z^{-k}. \lb{5.66} 
\end{align}
Comparing coefficients $z^{-k}$, $k\in\bbN$, then yields
\eqref{5.62a} and \eqref{5.62b}. In the special case $k=1$ one infers 
\begin{equation}
-\sum_{j=1}^n \cU_j(x)=-\frac{1}{2}\bigg(\sum_{\ell=0}^{2n}
E_\ell\bigg)\cI_m + \hatt\gR_{1,\Sigma}(x) \lb{5.67}
\end{equation}
and since $\hatt\gR_{1,\Sigma}=\cQ_\Sigma/2$ by \eqref{5.10b}, the trace
formula \eqref{5.63} for $\cQ_\Sigma$ results. \eqref{5.62c},
\eqref{5.62d}, and \eqref{5.64} are proved analogously. By
\eqref{5.61}, $\cU_j\in L^\infty(\bbR)^{m\times m}$, $1\leq j\leq n$ and
hence 
\begin{equation}
\hatt\gR_k\in C^\infty(\bbR)^{m\times m}\cap L^\infty (\bbR)^{m\times
m}, \quad k\in\bbN_0 
\lb{5.67a}
\end{equation}
(since $\cQ_\Sigma\in C^\infty(\bbR)^{m\times m}$ by Theorem
\ref{t4.9}\,(iii)). An analysis of the recursion relation \eqref{2.51}
for $\cM_{\pm,k}$ combined with \eqref{5.7}, \eqref{5.9} then proves
that $\hatt\gR_k$ is of the form  
\begin{equation}  
\hatt\gR_{k}=d_k \cQ^{(2k-2)} + \cR_k(\cQ^{(2k-4)},\dots), \quad k\geq 2,
\lb{5.67b}
\end{equation}
with $d_k\in\bbR$ appropriate constants and $\cR_k$ abbreviating a
differential polynomial in $\cQ$ which contains $\cQ^{(2k-4)}$ as the
highest derivative of $\cQ$. Hence one infers \eqref{5.65}.
\end{proof}
The factorizations \eqref{5.58}, \eqref{5.59}, eigenvalue distributions
\eqref{5.60}--\eqref{5.62}, and trace formulas
\eqref{5.62a}--\eqref{5.64} are extensions of well-known formulas in the
scalar case $m=1$ (see, e.g., \cite{Du75}, \cite{DMN76}, \cite{GHSZ95},
\cite{GRT96}, \cite{IM75}, \cite{LS88}).

Finally, the property $\cQ_\Sigma\in C^\infty(\bbR)^{m\times m}$ in
Theorem \ref{t4.9}\,(iii) can be improved upon by using the system
\eqref{4.58}, \eqref{4.59}, \eqref{4.60}, and \eqref{4.61}. In fact,
writing
\begin{align}
\cF_{n,\Sigma}(z,x)&=\sum_{\ell=0}^n \cF_{n-\ell,\Sigma}(x) z^\ell,
\quad \cF_{0,\Sigma}(x)=\cI_m, \lb{4.71} \\
\cG_{p,n-1,\Sigma}(z,x)&=\sum_{\ell=0}^{n-1} \cG_{p,n-1-\ell,\Sigma}(x)
z^\ell, \quad p=1,2, \lb{4.72} \\
\cH_{n+1,\Sigma}(z,x)&=\sum_{\ell=0}^{n+1} \cH_{n+1-\ell,\Sigma}(x)
z^\ell, \quad \cH_{0,\Sigma}(x)=\cI_m, \lb{4.73} 
\end{align}
one obtains the following result.

\begin{lemma} \lb{l5.8}
Assume Hypothesis \ref{h4.4} and let $(z,x)\in\bbC\times\bbR$. Then the
coefficients in \eqref{4.71}--\eqref{4.73} satisfy the autonomous
nonlinear first-order system
\begin{align}
\cF_{\ell,\Sigma}'&=-(\cG_{1,\ell-1,\Sigma}
+\cG_{2,\ell-1,\Sigma}), \quad 1\leq \ell\leq n, \lb{4.84} \\
\cG_{1,\ell,\Sigma}'&=-(\cF_{1,\Sigma}
-\cH_{1,\Sigma})\cF_{\ell+1,\Sigma}+\cF_{\ell+2,\Sigma}
-\cH_{\ell+2,\Sigma}, \quad 0\leq \ell\leq n-1, \lb{4.85} \\
\cG_{2,\ell,\Sigma}'&=-\cF_{\ell+1,\Sigma}(\cF_{1,\Sigma}
-\cH_{1,\Sigma})+\cF_{\ell+2,\Sigma}
-\cH_{\ell+2,\Sigma}, \quad 0\leq \ell\leq n-1, \lb{4.86} \\
\cH_{\ell,\Sigma}'&=\cG_{1,\ell-1,\Sigma}+\cG_{2,\ell-1,\Sigma}
-\cG_{1,\ell-2,\Sigma}(\cF_{1,\Sigma}-\cH_{1,\Sigma})
-(\cF_{1,\Sigma}-\cH_{1,\Sigma})\cG_{2,\ell-2,\Sigma}, \no \\
& \hspace*{8cm} 1\leq \ell\leq n+1, \lb{4.87} \\
\cF_{n+1,\Sigma}&=0, \quad \cG_{p,n,\Sigma}=\cG_{p,-1,\Sigma}=0, \quad
p=1,2. \lb{4.88} 
\end{align}
Moreover, $\cF_{\ell,\Sigma}$, $0\leq\ell\leq n$, $\cG_{p,\ell,\Sigma}$, 
$0\leq\ell\leq n-1$, $p=1,2$, $\cH_{\ell,\Sigma}$, $ 0\leq\ell\leq n+1$, 
and hence $\cQ_\Sigma$, are all analytic in an open neighborhood
containing the real axis.
\end{lemma}
\begin{proof}
Inserting \eqref{4.71}--\eqref{4.73} into \eqref{4.58},
\eqref{4.59}, \eqref{4.60}, and \eqref{4.61} yields
\begin{align}
\cF_{\ell,\Sigma}'&=-(\cG_{1,\ell-1,\Sigma}
+\cG_{2,\ell-1,\Sigma}), \quad 1\leq \ell\leq n, \lb{4.74} \\
\cG_{1,\ell,\Sigma}'&=-\cQ_\Sigma\cF_{\ell+1,\Sigma}+\cF_{\ell+2,\Sigma}
-\cH_{\ell+2,\Sigma}, \quad 0\leq \ell\leq n-1, \lb{4.75} \\
\cG_{2,\ell,\Sigma}'&=-\cF_{\ell+1,\Sigma}\cQ_\Sigma+\cF_{\ell+2,\Sigma}
-\cH_{\ell+2,\Sigma}, \quad 0\leq \ell\leq n-1, \lb{4.76} \\
\cH_{\ell,\Sigma}'&=\cG_{1,\ell-1,\Sigma}+\cG_{2,\ell-1,\Sigma}
-(\cG_{1,\ell-2,\Sigma}\cQ_\Sigma +\cQ_\Sigma\cG_{2,\ell-2,\Sigma}), 
\quad 1\leq \ell\leq n+1, \lb{4.77} \\
\cF_{n+1,\Sigma}&=0, \quad \cG_{p,n,\Sigma}=\cG_{p,-1,\Sigma}=0, \quad
p=1,2, \lb{4.78} \\
\cQ_\Sigma&=\cF_{1,\Sigma}-\cH_{1,\Sigma}. \lb{4.79}
\end{align}
Insertion of \eqref{4.79} into \eqref{4.75}--\eqref{4.77} yields the
autonomous nonlinear first-order system \eqref{4.84}--\eqref{4.88}. 
Given the initial conditions
\begin{align}
&\cF_{\ell,\Sigma}(x_0), \quad 1\leq \ell \leq n, \no \\
&\cG_{p,\ell,\Sigma}(x_0), \quad p=1,2, \;\, 0\leq \ell \leq n-1,
\lb{5.91}
\\ &\cH_{\ell,\Sigma}(x_0), \quad 1\leq \ell \leq n+1, \no 
\end{align}
determined by \eqref{4.71}--\eqref{4.73} and the half-line
Weyl--Titchmarsh matrices $\cM_\pm(z,x_0)$ in \eqref{4.40}, the maximal
interval of existence of the solution of the autonomous system
\eqref{4.84}--\eqref{4.88}, \eqref{5.91} is all of $[0,\infty)$ and
$(-\infty,0]$, and hence all of $\bbR$ (cf.\ \cite[p.\ 18]{Ha80}),
applying \eqref{5.65}. Thus, one recovers the
$C^\infty(\bbR)^{m\times m}$-property $\cF_\ell$, $\cG_{p,\ell}$, $p=1,2$,
$\cH_\ell$, and hence that of $\cQ_\Sigma$. Moreover, since the Picard
iterations are convergent in sufficiently small circles in $\bbC$
centered around each $x\in\bbR$ (cf.\ \cite[Sect.\ 2.3]{Hi97}), the unique
solution obtained by these Picard iterations is analytic in each of these
circles. 
\end{proof}

We conclude with a remark that puts the construction of $\cQ_\Sigma$ in 
Section \ref{s4} into proper perspective.

\begin{remark} \lb{r5.9}
The simplest examples of potentials $\cQ_\Sigma$ described in Theorem
\ref{t4.9} are of the type
\begin{equation}
\cQ_\Sigma(x)=\diag (q_{1,\Sigma}(x),\dots,q_{m,\Sigma}(x)), \lb{5.92}
\end{equation}
where $q_{j,\Sigma}$, $1\leq j\leq m$, are isospectral algebro-geometric
finite-band potentials associated with scalar Schr\"odinger
operators in $L^2(\bbR)$ with spectrum $\Sigma$. The next simplest and
closely related set of such examples for $\cQ_\Sigma$ then will be of the
type
\begin{equation}
\cQ_\Sigma(x)=\cU \diag (q_{1,\Sigma}(x),\dots,q_{m,\Sigma}(x))\cU^{-1},
\lb{5.93}
\end{equation}
where $\cU\in\bbC^{m\times m}$ is a unitary $m\times m$ matrix independent
of $x\in\bbR$. At first sight one might think that perhaps all matrix
potentials $\cQ_\Sigma$ are of the form \eqref{5.93}. That this is
certainly not the case will be argued next. Indeed, assuming that
\eqref{5.93} holds for some unitary $m\times m$ matrix $\cU$ independent
of $x$, one infers that 
\begin{equation}
\cQ_\Sigma^{(r)}(x)=\cU \diag
(q_{1,\Sigma}^{(r)}(x),\dots,q_{m,\Sigma}^{(r)}(x))\cU^{-1} 
\text{ for all $r\in\bbN_0$}. \lb{5.94}
\end{equation}
Since by Remark \ref{r2.11}\,(i) and Lemma \ref{l5.4} all coefficients of
$\cF_{n,\Sigma}$, $\cG_{p,n-1,\Sigma}$, $p=1,2$, and $\cH_{n+1,\Sigma}$
are differential polynomials with respect to $\cQ_\Sigma$, \eqref{5.94}
implies
\begin{align}
\cU\cF_{n,\Sigma}(z,x)\cU^{-1}&=\diag(F_{n,\Sigma,1}(z,x),\dots,
F_{n,\Sigma,m}(z,x)), \lb{5.95} \\
\cU\cG_{p,n-1,\Sigma}(z,x)\cU^{-1}&=\diag(G_{p,n-1,\Sigma,1}(z,x),\dots,
G_{p,n-1,\Sigma,m}(z,x)), \quad p=1,2, \lb{5.96} \\
\cU\cH_{n+1,\Sigma}(z,x)\cU^{-1}&=\diag(H_{n+1,\Sigma,1}(z,x),\dots,
H_{n+1,\Sigma,m}(z,x)), \lb{5.97} \\
\cU\cM_{\pm,\Sigma}(z,x)\cU^{-1}&=\diag(m_{\pm,\Sigma,1}(z,x),\dots,
m_{\pm,\Sigma,m}(z,x)). \lb{5.98} 
\end{align}
Consequently, one obtains for all $z,z'\in\bbC_+$, $x,x'\in\bbR$,
\begin{align}
[\cQ_\Sigma^{(r)} (x),\cQ_\Sigma^{(r)} (x')]&=0, \quad
r\in\bbN_0, \lb{5.98a}
\\ [\cF_{n,\Sigma}(z,x),\cF_{n,\Sigma}(z',x')]&=0, \lb{5.99} \\
[\cG_{p,n-1,\Sigma}(z,x),\cG_{p,n-1,\Sigma}(z',x')]&=0, \quad
p=1,2, \lb{5.100} \\ 
[\cH_{n+1,\Sigma}(z,x),\cH_{n+1,\Sigma}(z',x')]&=0, \lb{5.101} \\
[\cM_{\pm,\Sigma}(z,x),\cM_{\pm,\Sigma}(z',x')]&=0. \lb{5.102} 
\end{align}
In particular,
\begin{equation}
[\cF_{n,\Sigma}(z,x_0),\cF_{n,\Sigma}(z',x_0)]=0 
\text{ for all $z,z'\in\bbC_+$.} \lb{5.103} 
\end{equation}
Thus, whenever
\begin{equation}
[\cF_{n,\Sigma}(z,x_0),\cF_{n,\Sigma}(z',x_0)]\neq 0 
\text{ for some $z,z'\in\bbC_+$} \lb{5.104} 
\end{equation}
$($which can easily be arranged for $n\geq 2$$)$, \eqref{5.93} cannot hold
for a unitary $m\times m$ matrix $\cU$ independent of $x$.
\end{remark}

Additional results, including extensions of Borg's and Hochstadt's 
theorems in the special cases $n=0,1$, respectively, will appear
in \cite{BGMS02}.

\bigskip
\noindent {\bf Acknowledgements.} We thank Robert Carlson, Konstantin
Makarov, Mark Malamud, Fedor Rofe-Beketov, and  Barry Simon for many
helpful discussions on the material presented in this paper. In 
particular, we are grateful to Konstantin Makarov for a critical reading
of  this manuscript.


\end{document}